\documentclass[a4paper,twoside]{article}
\usepackage{amsfonts}
\usepackage{amsmath}
\usepackage{amsthm}
\let\topsmash\smash
\numberwithin{equation}{section}
\renewcommand{\[}{\lfloor}
\renewcommand{\]}{\rfloor}
\newcommand{\<}{\langle}
\renewcommand{\>}{\rangle}
\renewcommand{\kappa}{\varkappa}
\renewcommand{\epsilon}{\varepsilon}
\renewcommand{\phi}{\varphi}
\newcommand{\vt}{\vartheta}
\renewcommand{\d}{\mathrm{d}}
\newcommand{\wt}{\widetilde}

\DeclareMathOperator{\ord}{ord}
\DeclareMathOperator{\Real}{Re}
\renewcommand{\Re}{\Real}
\DeclareMathOperator{\Imag}{Im}
\renewcommand{\Im}{\Imag}
\DeclareMathOperator{\Res}{Res}
\DeclareMathOperator{\id}{id}
\newcommand{\rom}[1]{{\rm#1}}
\newcommand{\doublesb}[2]{_{\genfrac{}{}{0pt}{1}{#1}{#2}}}
\newcommand{\ba}{\boldsymbol a}
\newcommand{\bb}{\boldsymbol b}
\newcommand{\bc}{\boldsymbol c}
\newcommand{\be}{\boldsymbol e}
\newcommand{\bh}{\boldsymbol h}
\newcommand{\balpha}{\boldsymbol\alpha}
\newcommand{\bbeta}{\boldsymbol\beta}
\newcommand{\Beta}{\boldsymbol\eta}
\newcommand{\fa}{\mathfrak a}
\newcommand{\fb}{\mathfrak b}
\newcommand{\fh}{\mathfrak h}
\newcommand{\fq}{\mathfrak q}
\newcommand{\fA}{\mathfrak A}
\newcommand{\fS}{\mathfrak S}
\newcommand{\fG}{\mathfrak G}
\newcommand{\sL}{\mathcal L}
\newcommand{\sM}{\mathcal M}
\begin{document}
\newtheorem{theorem}{Theorem}
\newtheorem{lemma}{Lemma}
\newtheorem{proposition}{Proposition}
\theoremstyle{remark}
\newtheorem*{rem}{Remark}
\newtheorem*{observation}{Observation}
\theoremstyle{plain}
\newtheorem*{conj}{Conjecture}
\title{Arithmetic of linear forms\\ involving odd zeta values%
}
\author{Wadim Zudilin%
\thanks{Moscow Lomonosov State University,
Department of Mechanics and Mathematics,
Vorobiovy Gory, GSP-2, 119992 Moscow, Russia.
e-mail: \texttt{wadim@ips.ras.ru}}}
\date{E-print \texttt{math.NT/0206176}\\[1.5mm] August 2001}
\maketitle
\pagestyle{myheadings}
\markboth{W.~Zudilin}{Arithmetic of linear forms involving odd zeta values}
\thispagestyle{empty}

\begin{abstract}
A general hypergeometric construction of linear forms in (odd) zeta values
is presented. The construction allows to recover the records of Rhin and
Viola for the irrationality measures of $\zeta(2)$ and $\zeta(3)$,
as well as to explain Rivoal's recent result (\texttt{math.NT/0008051})
on infiniteness of irrational numbers in the set of odd zeta values,
and to prove that at least one of the four numbers
$\zeta(5)$, $\zeta(7)$, $\zeta(9)$, and $\zeta(11)$ is irrational.
\end{abstract}


\section{Introduction}
\label{sec:1}
The story exposed in this paper starts in 1978, when R.~Ap\'ery~\cite{Ap}
gave a surprising sequence of exercises demonstrating the
irrationality of~$\zeta(2)$ and~$\zeta(3)$. (For a nice explanation
of Ap\'ery's discovery we refer to the review~\cite{Po}.)
Although the irrationality of the even zeta values $\zeta(2),\zeta(4),\dots$
for that moment was a classical result (due to L.~Euler and F.~Lindemann),
Ap\'ery's proof allows one to obtain a {\it quantitative\/} version
of his result, that is, to evaluate irrationality exponents:
\begin{equation}
\mu(\zeta(2))\le11.85078\dots,
\qquad
\mu(\zeta(3))\le13.41782\dots\,.
\label{eq:1.1}
\end{equation}
As usual, a value $\mu=\mu(\alpha)$ is said to be the {\it irrationality
exponent\/} of an irrational number~$\alpha$ if $\mu$~is the least
possible exponent such that for any $\epsilon>0$ the inequality
$$
\biggl|\alpha-\frac pq\biggr|\le\frac1{q^{\mu+\epsilon}}
$$
has only finitely many solutions in integers $p$ and~$q$ with $q>0$.
The estimates~\eqref{eq:1.1} `immediately' follow from the asymptotics
of Ap\'ery's rational approximations to~$\zeta(2)$ and~$\zeta(3)$,
and the original method of evaluating the asymptotics is
based on second order difference equations with polynomial
coefficients, with Ap\'ery's approximants as their solutions.

A few months later,
F.~Beukers~\cite{Be} interpretated Ap\'ery's sequence
of rational approximations to $\zeta(2)$ and $\zeta(3)$
in terms of multiple integrals and Legendre polynomials.
This approach was continued in later
works~\cite{DV,Ru}, \cite{Ha1}--\cite{Ha5},
\cite{HMV}, \cite{RV1}--\cite{RV3}
and yielded some new evaluations of the irrationality exponents
for $\zeta(2)$, $\zeta(3)$, and other mathematical constants.
Improvements of irrationality measures (i.e., upper bounds for
irrationality exponents) for mathematical constants
are closely related to another {\it arithmetic\/} approach, of eliminating
extra prime numbers in binomials, introduced after G.\,V.~Chudnovsky~\cite{Ch}
by E.\,A.~Rukhadze~\cite{Ru} and studied in detail by M.~Hata~\cite{Ha1}.
For example, the best known estimate
for the irrationality exponent of~$\log2$
(this constant sometimes is regarded as a convergent
analogue of~$\zeta(1)$\,) stated by Rukhadze~\cite{Ru} in~1987 is
\begin{equation}
\mu(\log2)\le3.891399\dots;
\label{eq:1.2}
\end{equation}
see also~\cite{Ha1} for the explicit value of the constant
on the right-hand side of~\eqref{eq:1.2}.
A further generalization of both
the multiple integral approach and the arithmetic approach
brings one to the group structures of G.~Rhin and
C.~Viola~\cite{RV2,RV3}; their method yields
the best known estimates for the irrationality exponents
of $\zeta(2)$ and~$\zeta(3)$:
\begin{equation}
\mu(\zeta(2))\le5.441242\dots,
\qquad
\mu(\zeta(3))\le5.513890\dots,
\label{eq:1.3}
\end{equation}
and gives another interpretation~\cite{Vi}
of Rukhadze's estimate~\eqref{eq:1.2}.

On the other hand, Ap\'ery's phenomenon was interpretated
by L.\,A.~Gutnik~\cite{Gu} in terms of complex contour integrals,
i.e., Meijer's $G$-functions. This approach allowed the author of~\cite{Gu}
to prove several partial results on the irrationality of certain quantities
involving $\zeta(2)$ and~$\zeta(3)$. By the way of a study
of Gutnik's approach, Yu.\,V.~Nesterenko~\cite{Ne1} proposed
a new proof of Ap\'ery's theorem and discovered
a new continuous fraction expansion for~$\zeta(3)$. In~\cite{FN}, p.~126,
a problem of finding an `elementary' proof
of the irrationality of~$\zeta(3)$ is stated since evaluating asymptotics
of multiple integrals via the Laplace method in~\cite{Be} or
complex contour integrals via the saddle-point method in~\cite{Ne1}
is far from being simple. Trying to solve this problem, K.~Ball
puts forward a well-poised hypergeometric series, which produces linear
forms in~$1$ and~$\zeta(3)$ only and can be evaluated
by elementary means;
however, its `obvious' arithmetic does not allow one to prove
the irrationality of~$\zeta(3)$. T.~Rivoal~\cite{Ri1} has realized
how to generalize Ball's linear form in the spirit of
Nikishin's work~\cite{Ni} and to use well-poised
hypergeometric series in the study of the irrationality of odd
zeta values $\zeta(3),\zeta(5),\dots$; in particular, he is able
to prove~\cite{Ri1} that there are infinitely many irrational numbers
in the set of the odd zeta values. A further generalization of the method
in the spirit of~\cite{Gu,Ne1} via the use of well-poised
Meijer's $G$-functions allows Rivoal~\cite{Ri4} to demonstrate
the irrationality of at least one of the nine numbers
$\zeta(5),\zeta(7),\dots,\zeta(21)$. Finally, this
author~\cite{Zu1}--\cite{Zu4} refines the results of
Rivoal~\cite{Ri1}--\cite{Ri4} by an application of
the arithmetic approach.

Thus, one can recognise (at least) two different languages used for
an explanation why $\zeta(3)$~is irrational, namely, multiple integrals
and complex contour integrals (or series of hypergeometric type).
Both languages lead us to quantitative and qualitative results
on the irrationality of zeta values and other mathematical constants,
and it would be nice to form a dictionary for translating terms from
one language into another. An approach to such a translation
has been recently proposed by Nesterenko~\cite{Ne2,Ne3}.
He has proved a general theorem that expresses contour integrals
in terms of multiple integrals, and vice versa. He also
suggests a method of constructing linear forms in values
of polylogarithms (and, as a consequence, linear forms in zeta values)
that generalizes the language
of~\cite{Ni,Gu,Ne1} and, on the other hand,
of~\cite{Be}, \cite{Ha1}--\cite{Ha5}, \cite{RV1}--\cite{RV3}
and takes into account both arithmetic and analytic evaluations
of the corresponding linear forms.

The aim of this paper is to explain the group structures used
for evaluating the irrationality exponents~\eqref{eq:1.2},~\eqref{eq:1.3}
via Nesterenko's method, as well as to present a new result
on the irrationality of the odd zeta values inspired by Rivoal's
construction and possible generalizations of the Rhin--Viola approach.
This paper is organized as follows. In Sections~\ref{sec:2}--\ref{sec:5}
we explain
in details the group structure of Rhin and Viola for~$\zeta(3)$;
we do not use Beukers' type integrals as in~\cite{RV3} for this,
but with the use of Nesterenko's theorem we explain
all stages of our construction in terms of their doubles from~\cite{RV3}.
Section~\ref{sec:6} gives a brief overview of the group structure
for~$\zeta(2)$ from~\cite{RV2}.
Section~\ref{sec:7} is devoted to a study of the arithmetic of rational
functions appearing naturally as `bricks' of general Nesterenko's
construction~\cite{Ne3}. In Section~\ref{sec:8} we explain
the well-poised hypergeometric origin of
Rivoal's construction and improve the previous result
from~\cite{Ri4,Zu4}
on the irrationality of $\zeta(5),\zeta(7),\dots$; namely,
we state that at least one of the four numbers
$$
\zeta(5), \; \zeta(7), \; \zeta(9), \; \mbox{and} \; \zeta(11)
$$
is irrational. Although the success of our new result from Section~\ref{sec:8}
is due to the arithmetic approach, in Section~\ref{sec:9} we present
possible group structures for linear forms in~$1$ and odd zeta values;
these groups may become useful, provided that some arithmetic condition
(which we indicate explicitly) holds.

\medskip
This work would be not possible without a permanent attention of
Professor Yu.\,V.~Nesterenko. I would like to express my deep
gratitude to him. I am thankful to T.~Rivoal for giving
me the possibility to look through his Ph.~D.~thesis~\cite{Ri3},
which contains a lot of fruitful ideas exploited in this work.

This research was carried out with
the partial support of the INTAS--RFBR grant no.~IR-97-1904.

\section{Analytic construction\\ of linear forms in $1$ and $\zeta(3)$}
\label{sec:2}
Fix a set of integral parameters
\begin{equation}
(\ba,\bb)
=\biggl(\begin{array}{rrrr}
a_1, & a_2, & a_3, & a_4 \\
b_1, & b_2, & b_3, & b_4
\end{array}\biggr)
\label{eq:2.1}
\end{equation}
satisfying the conditions
\begin{gather}
\{b_1,b_2\}\le\{a_1,a_2,a_3,a_4\}<\{b_3,b_4\},
\label{eq:2.2}
\\
a_1+a_2+a_3+a_4\le b_1+b_2+b_3+b_4-2,
\label{eq:2.3}
\end{gather}
and consider the rational function
\begin{equation}
\begin{split}
R(t)
=R(\ba,\bb;t)
&:=\frac{(b_3-a_3-1)!\,(b_4-a_4-1)!}{(a_1-b_1)!\,(a_2-b_2)!}
\\ &\phantom:\qquad\times
\frac{\Gamma(t+a_1)\,\Gamma(t+a_2)\,\Gamma(t+a_3)\,\Gamma(t+a_4)}
{\Gamma(t+b_1)\,\Gamma(t+b_2)\,\Gamma(t+b_3)\,\Gamma(t+b_4)}
\\ &\phantom:
=\prod_{j=1}^4R_j(t),
\end{split}
\label{eq:2.4}
\end{equation}
where
\begin{equation}
R_j(t)=\begin{cases}
\dfrac{(t+b_j)(t+b_j+1)\dotsb(t+a_j-1)}{(a_j-b_j)!}
& \mbox{if $a_j\ge b_j$ (i.e., $j=1,2$)}, \\
\dfrac{(b_j-a_j-1)!}{(t+a_j)(t+a_j+1)\dotsb(t+b_j-1)}
& \mbox{if $a_j<b_j$ (i.e., $j=3,4$)}.
\end{cases}
\label{eq:2.5}
\end{equation}
By condition~\eqref{eq:2.3} we obtain
\begin{equation}
R(t)=O(t^{-2})
\qquad\mbox{as}\quad t\to\infty;
\label{eq:2.6}
\end{equation}
moreover, the function $R(t)$ has zeros of the second order
at the integral points~$t$ in the interval
$$
-\min\{a_1,a_2,a_3,a_4\}<t\le-\max\{b_1,b_2\}.
$$
Therefore, the numerical series $\sum_{t=t_0}^\infty R'(t)$
with $t_0=1-\max\{b_1,b_2\}$ converges absolutely,
and the quantity
\begin{equation}
G(\ba,\bb)
:=-(-1)^{b_1+b_2}\sum_{t=t_0}^\infty R'(t)
\label{eq:2.7}
\end{equation}
is well-defined; moreover, we can start the summation
on the right-hand side of~\eqref{eq:2.7} from any integer~$t_0$
in the interval
\begin{equation}
1-\min\{a_1,a_2,a_3,a_4\}\le t_0\le 1-\max\{b_1,b_2\}.
\label{eq:2.8}
\end{equation}
The number~\eqref{eq:2.7} is a linear form in~$1$ and~$\zeta(3)$
(see Lemma~\ref{lem:4} below), and we devote the rest of this section
to a study of the arithmetic (i.e., the denominators of the coefficients)
of this linear form.

To the data~\eqref{eq:2.1} we assign the ordered set $(\ba^*,\bb^*)$; namely,
\begin{equation}
\begin{gathered}
\{b_1^*,b_2^*\}=\{b_1,b_2\},
\quad
\{a_1^*,a_2^*,a_3^*,a_4^*\}=\{a_1,a_2,a_3,a_4\},
\\
\{b_3^*,b_4^*\}=\{b_3,b_4\},
\qquad
b_1^*\le b_2^*
\le a_1^*\le a_2^*\le a_3^*\le a_4^*
<b_3^*\le b_4^*,
\end{gathered}
\label{eq:2.9}
\end{equation}
hence the interval~\eqref{eq:2.8} for~$t_0$ can be written as follows:
$$
1-a_1^*\le t_0\le 1-b_2^*.
$$
By~$D_N$ we denote the least common multiple of numbers $1,2,\dots,N$.

\begin{lemma}
\label{lem:1}
For $j=1,2$ there hold the inclusions
\begin{equation}
R_j(t)\big|_{t=-k}\in\mathbb Z,
\quad
D_{a_j-b_j}\cdot R_j'(t)\big|_{t=-k}\in\mathbb Z,
\qquad k\in\mathbb Z.
\label{eq:2.10}
\end{equation}
\end{lemma}

\begin{proof}
The inclusions~\eqref{eq:2.10} immediately follow from the well-known
properties of the {\it integral-valued polynomials\/}
(see, e.g., \cite{Zu5}, Lemma~7),
which are $R_1(t)$ and~$R_2(t)$.
\end{proof}

The analogue of Lemma~\ref{lem:1}
for rational functions $R_3(t),R_4(t)$ from~\eqref{eq:2.5}
is based on the following assertion
combining the arithmetic schemes of Ni\-ki\-shin~\cite{Ni}
and Rivoal~\cite{Ri1}.

\begin{lemma}[\textrm{\rm\cite{Zu3}, Lemma~1.2}]
\label{lem:2}
Assume that for some polynomial~$P(t)$ of degree not greater than~$n$
the rational function
$$
Q(t)=\frac{P(t)}{(t+s)(t+s+1)\dotsb(t+s+n)}
$$
\rom(in a not necesarily uncancellable presentation\rom)
satisfies the conditions
$$
Q(t)(t+k)\big|_{t=-k}\in\mathbb Z,
\qquad k=s,s+1,\dots,s+n.
$$
Then for all non-negative integers~$l$ there hold the inclusions
$$
\frac{D_n^l}{l!}\cdot\bigl(Q(t)(t+k)\bigr)^{(j)}\big|_{t=-k}\in\mathbb Z,
\qquad k=s,s+1,\dots,s+n.
$$
\end{lemma}

\begin{lemma}
\label{lem:3}
For $j=3,4$ there hold the inclusions
\begin{gather}
\bigl(R_j(t)(t+k)\bigr)\big|_{t=-k}\in\mathbb Z,
\qquad k\in\mathbb Z,
\label{eq:2.11}
\\
\begin{gathered}
D_{b_4^*-\min\{a_j,a_3^*\}-1}
\cdot\bigl(R_j(t)(t+k)\bigr)'\big|_{t=-k}\in\mathbb Z,
\\
k\in\mathbb Z, \quad k=a_3^*,a_3^*+1,\dots,b_4^*-1.
\end{gathered}
\label{eq:2.12}
\end{gather}
\end{lemma}

\begin{proof}
The inclusions~\eqref{eq:2.11} can be verified by direct calculations:
$$
\bigl(R_j(t)(t+k)\bigr)\big|_{t=-k}=\begin{cases}
(-1)^{k-a_j}\dfrac{(b_j-a_j-1)!}{(k-a_j)!\,(b_j-k-1)!}
\\ \phantom0\quad
\mbox{if $k=a_j,a_j+1,\dots,b_j-1$},
\\
0 \quad\mbox{otherwise}.
\end{cases}
$$
To prove the inclusions~\eqref{eq:2.12}
we apply Lemma~\ref{lem:2} with $l=1$
to the function~$R_j(t)$ multiplying its numerator and denominator
if necesary by the factor
$(t+a_3^*)\dotsb\linebreak[2](t+a_j-1)$
if $a_j>a_3^*$ and by
$(t+b_j)\dotsb(t+b_4^*-1)$ if $b_j<b_4^*$.
\end{proof}

\begin{lemma}
\label{lem:4}
The quantity~\eqref{eq:2.7} is a linear form in~$1$ and~$\zeta(3)$
with rational coefficients\rom:
\begin{equation}
G(\ba,\bb)=2A\zeta(3)-B;
\label{eq:2.13}
\end{equation}
in addition\rom,
\begin{equation}
A\in\mathbb Z, \qquad
D_{b_4^*-a_1^*-1}^2
\cdot D_{\max\{a_1-b_1,a_2-b_2,b_4^*-a_3-1,b_4^*-a_4-1,b_3^*-a_1^*-1\}}
\cdot B\in\mathbb Z.
\label{eq:2.14}
\end{equation}
\end{lemma}

\begin{proof}
The rational function~\eqref{eq:2.4} has poles
at the points $t=-k$,
where $k=a_3^*,a_3^*+1,\dots,b_4^*-1$;
moreover, the points $t=-k$,
where $k=a_4^*,a_4^*+1,\dots,b_3^*-1$,
are poles of the second order.
Hence the expansion of the rational function~\eqref{eq:2.4}
in a sum of partial fractions has the form
\begin{equation}
R(t)=\sum_{k=a_4^*}^{b_3^*-1}\frac{A_k}{(t+k)^2}
+\sum_{k=a_3^*}^{b_4^*-1}\frac{B_k}{t+k},
\label{eq:2.15}
\end{equation}
where the coefficients $A_k$ and~$B_k$ in~\eqref{eq:2.15}
can be calculated by the formulae
$$
\begin{alignedat}{2}{2}
A_k&=\bigl(R(t)(t+k)^2\bigr)\big|_{t=-k}, \qquad&
k&=a_4^*,a_4^*+1,\dots,b_3^*-1,
\\
B_k&=\bigl(R(t)(t+k)^2\bigr)'\big|_{t=-k}, \qquad&
k&=a_3^*,a_3^*+1,\dots,b_4^*-1.
\end{alignedat}
$$
Expressing the function $R(t)(t+k)^2$ as
$$
R_1(t)\cdot R_2(t)\cdot R_3(t)(t+k)\cdot R_4(t)(t+k)
$$
for each~$k$ and applying the Leibniz rule for differentiating
a product, by Lemmas~\ref{lem:1} and~\ref{lem:3} we obtain
\begin{equation}
\begin{alignedat}{2}{2}
A_k&\in\mathbb Z,
\qquad&
k&=a_4^*,a_4^*+1,\dots,b_3^*-1,
\\
D_{\max\{a_1-b_1,a_2-b_2,b_4^*-a_3-1,b_4^*-a_4-1\}}
\cdot B_k&\in\mathbb Z,
\qquad&
k&=a_3^*,a_3^*+1,\dots,b_4^*-1
\end{alignedat}
\label{eq:2.16}
\end{equation}
(where we use the fact that $\min\{a_j,a_3^*\}\le a_j$
for at least one $j\in\{3,4\}$).

By~\eqref{eq:2.6} there holds
$$
\sum_{k=a_3^*}^{b_4^*-1}B_k
=\sum_{k=a_3^*}^{b_4^*-1}\Res_{t=-k}R(t)
=-\Res_{t=\infty}R(t)=0.
$$
Hence, setting $t_0=1-a_1^*$ in~\eqref{eq:2.7}
and using the expansion~\eqref{eq:2.15} we obtain
\begin{align*}
(-1)^{b_1+b_2}G(\ba,\bb)
&=\sum_{t=1-a_1^*}^\infty
\biggl(\sum_{k=a_4^*}^{b_3^*-1}\frac{2A_k}{(t+k)^3}
+\sum_{k=a_3^*}^{b_4^*-1}\frac{B_k}{(t+k)^2}\biggr)
\\
&=2\sum_{k=a_4^*}^{b_3^*-1}
A_k\biggl(\sum_{l=1}^\infty-\sum_{l=1}^{k-a_1^*}\biggr)\frac1{l^3}
+\sum_{k=a_3^*}^{b_4^*-1}
B_k\biggl(\sum_{l=1}^\infty-\sum_{l=1}^{k-a_1^*}\biggr)\frac1{l^2}
\\
&=2\sum_{k=a_4^*}^{b_3^*-1}A_k\cdot\zeta(3)
-\biggl(2\sum_{k=a_4^*}^{b_3^*-1}A_k\sum_{l=1}^{k-a_1^*}\frac1{l^3}
+\sum_{k=a_3^*}^{b_4^*-1}B_k\sum_{l=1}^{k-a_1^*}\frac1{l^2}\biggr)
\\
&=2A\zeta(3)-B.
\end{align*}
The inclusions~\eqref{eq:2.14} now follow from~\eqref{eq:2.16}
and the definition of the least common multiple:
\begin{alignat*}{2}
D_{b_4^*-a_1^*-1}^2\cdot\frac1{l^2}&\in\mathbb Z
\qquad&\mbox{for}\quad l&=1,2,\dots,b_4^*-a_1^*-1,
\\
D_{b_4^*-a_1^*-1}^2\cdot D_{b_3^*-a_1^*-1}
\cdot\frac1{l^3}&\in\mathbb Z
\qquad&\mbox{for}\quad l&=1,2,\dots,b_3^*-a_1^*-1.
\end{alignat*}
The proof is complete.
\end{proof}

Taking $a_1=a_2=a_3=a_4=1+n$, $b_1=b_2=1$, and $b_3=b_4=2+2n$
we obtain the original Ap\'ery's sequence
\begin{equation}
2A_n\zeta(3)-B_n
=-\sum_{t=1}^\infty\frac{\d}{\d t}
\biggl(\frac{(t-1)(t-2)\dotsb(t-n)}{t(t+1)\dotsb(t+n)}\biggr)^2,
\qquad n=1,2,\dots,
\label{eq:2.17}
\end{equation}
of rational approximations to~$\zeta(3)$
(cf.~\cite{Gu,Ne1}); Lemma~\ref{lem:4} implies that
$A_n\in\mathbb Z$ and $D_n^3\cdot B_n\in\mathbb Z$ in Ap\'ery's case.

\section{Integral presentations}
\label{sec:3}
The aim of this section is to prove two presentations
of the linear form~\eqref{eq:2.7}, \eqref{eq:2.13}:
as a complex contour integral
(in the spirit of~\cite{Gu,Ne1})
and as a real multiple integral
(in the spirit of~\cite{Be,Ha5,RV3}).

Consider another normalization of the rational function~\eqref{eq:2.4};
namely,
\begin{equation}
\wt R(t)
=\wt R(\ba,\bb;t)
:=\frac{\Gamma(t+a_1)\,\Gamma(t+a_2)\,\Gamma(t+a_3)\,\Gamma(t+a_4)}
{\Gamma(t+b_1)\,\Gamma(t+b_2)\,\Gamma(t+b_3)\,\Gamma(t+b_4)}
\label{eq:3.1}
\end{equation}
and the corresponding sum
\begin{equation}
\wt G(\ba,\bb)
:=-(-1)^{b_1+b_2}\sum_{t=t_0}^\infty\wt R'(t)
=\frac{(a_1-b_1)!\,(a_2-b_2)!}
{(b_3-a_3-1)!\,(b_4-a_4-1)!}G(\ba,\bb).
\label{eq:3.2}
\end{equation}
Note that the function~\eqref{eq:3.1} and the quantity~\eqref{eq:3.2}
do not depend on the order of numbers in the sets
$\{a_1,a_2,a_3,a_4\}$, $\{b_1,b_2\}$, and $\{b_3,b_4\}$, i.e.,
$$
\wt R(\ba,\bb;t)\equiv\wt R(\ba^*,\bb^*;t),
\qquad
\wt G(\ba,\bb)\equiv\wt G(\ba^*,\bb^*).
$$

\begin{lemma}
\label{lem:5}
There holds the formula
\begin{align}
\wt G(\ba,\bb)
&=\frac1{2\pi i}\int_{\sL}
\frac{\begin{aligned}
\Gamma(t+a_1)\,\Gamma(t+a_2)\,\Gamma(t+a_3)\,\Gamma(t+a_4)\,
\qquad \\[-3pt] \times
\Gamma(1-t-b_1)\,\Gamma(1-t-b_2)
\end{aligned}}
{\Gamma(t+b_3)\,\Gamma(t+b_4)}\,\d t
\nonumber\\
&=:G_{4,4}^{2,4}\biggl(1\biggm|\begin{array}{rrrr}
1-a_1, & 1-a_2, & 1-a_3, & 1-a_4 \\
1-b_1, & 1-b_2, & 1-b_3, & 1-b_4
\end{array}\biggr),
\label{eq:3.3}
\end{align}
where $\sL$~is a vertical line $\Re t=t_1$\rom, $1-a_1^*<t_1<1-b_2^*$\rom,
oriented from the bottom to the top\rom, and $G_{4,4}^{2,4}$~is Meijer's
$G$-function \rm(see~\cite{Lu}, Section~5.3).
\end{lemma}

\begin{proof}
The standard arguments (see, e.g., \cite{Gu}, \cite{Ne1}, Lemma~2,
or \cite{Zu3}, Lemma 2.4) show that the quantity~\eqref{eq:3.2}
presents the sum of the residues at the poles
$t=-b_2^*+1,-b_2^*+2,\dots$ of the function
\begin{align*}
&
-(-1)^{b_1+b_2}
\biggl(\frac\pi{\sin\pi t}\biggr)^2\wt R(t)
\\ &\qquad
=-(-1)^{b_1+b_2}
\biggl(\frac\pi{\sin\pi t}\biggr)^2
\frac{\Gamma(t+a_1)\,\Gamma(t+a_2)\,\Gamma(t+a_3)\,\Gamma(t+a_4)}
{\Gamma(t+b_1)\,\Gamma(t+b_2)\,\Gamma(t+b_3)\,\Gamma(t+b_4)}.
\end{align*}
It remains to observe that
\begin{equation}
\Gamma(t+b_j)\Gamma(1-t-b_j)=(-1)^{b_j}\frac\pi{\sin\pi t},
\qquad j=1,2,
\label{eq:3.4}
\end{equation}
and to identify the integral in~\eqref{eq:3.3} with Meijer's $G$-function.
This establishes formula~\eqref{eq:3.3}.
\end{proof}

The next assertion allows one to express the complex integral~\eqref{eq:3.3}
as a real multiple integral.

\begin{proposition}[\textrm{\rm Nesterenko's theorem \cite{Ne3}}]
\label{prop:1}
Suppose that $m\ge1$ and $r\ge0$~are integers\rom, $r\le m$\rom,
and that complex parameters $a_0,a_1,\dots,a_m$\rom, $b_1,\dots,b_m$
and a real number~$t_1<0$ satisfy the conditions
$$
\begin{gathered}
\Re b_k>\Re a_k>0, \qquad k=1,\dots,m,
\\
-\min_{0\le k\le m}\Re a_k<t_1
<\min_{1\le k\le r}\Re(b_k-a_k-a_0).
\end{gathered}
$$
Then for any $z\in\mathbb C\setminus(-\infty,0]$
there holds the identity
\begin{align*}
&
\idotsint\limits_{[0,1]^m}
\frac{\prod_{k=1}^mx_k^{a_k-1}(1-x_k)^{b_k-a_k-1}}
{\bigl((1-x_1)(1-x_2)\dotsb(1-x_r)+zx_1x_2\dotsb x_m\bigr)^{a_0}}
\,\d x_1\,\d x_2\dotsb\d x_m
\\ &\quad
=\frac{\prod_{k=r+1}^m\Gamma(b_k-a_k)}
{\Gamma(a_0)\cdot\prod_{k=1}^r\Gamma(b_k-a_0)}
\\ &\quad\quad\times
\frac1{2\pi i}\int_{t_1-i\infty}^{t_1+i\infty}
\frac{\prod_{k=0}^m\Gamma(a_k+t)
\cdot\prod_{k=1}^r\Gamma(b_k-a_k-a_0-t)}
{\prod_{k=r+1}^m\Gamma(b_k+t)}\,
\Gamma(-t)\,z^t\,\d t,
\end{align*}
where both integrals converge.
Here $z^t=e^{t\log z}$ and the logarithm takes
real values for real $z\in(0,+\infty)$.
\end{proposition}

We now recall that the family of linear forms
in~$1$ and $\zeta(3)$ considered in paper~\cite{RV3}
has the form
\begin{equation}
I(h,j,k,l,m,q,r,s)
=\iiint\limits_{[0,1]^3}
\frac{x^h(1-x)^ly^k(1-y)^sz^j(1-z)^q}{(1-(1-xy)z)^{q+h-r}}\,
\frac{\d x\,\d y\,\d z}{1-(1-xy)z}
\label{eq:3.5}
\end{equation}
and depends on eight non-negative integral parameters
connected by the additional conditions
\begin{equation}
h+m=k+r, \qquad j+q=l+s,
\label{eq:3.6}
\end{equation}
where the first condition in~\eqref{eq:3.6} determines the parameter~$m$
(which does not appear on the right-hand side of~\eqref{eq:3.5} explicitly),
while the second condition enables one to apply a complicated integral
transform~$\vt$, which rearranges all eight parameters.

\begin{lemma}
\label{lem:6}
The quantity~\eqref{eq:2.7} has the integral presentation
\begin{equation}
G(\ba,\bb)=I(h,j,k,l,m,q,r,s),
\label{eq:3.7}
\end{equation}
where the multiple integral on the right-hand side of~\eqref{eq:3.7}
is given by formula~\eqref{eq:3.5} and
\begin{equation}
\begin{alignedat}{4}{4}
h&=a_3-b_1, \quad&
j&=a_2-b_1, \quad&
k&=a_4-b_1, \quad&
l&=b_3-a_3-1,
\\
m&=a_4-b_2, \quad&
q&=a_1-b_2, \quad&
r&=a_3-b_2, \quad&
s&=b_4-a_4-1.
\end{alignedat}
\label{eq:3.8}
\end{equation}
\end{lemma}

\begin{proof}
By the change of variables $t\mapsto t-b_1+1$ in the complex
integral~\eqref{eq:3.3} and the application of Proposition~\ref{prop:1}
with $m=3$, $r=1$, and $z=1$ we obtain
\begin{align*}
\wt G(\ba,\bb)
&=\frac{(a_1-b_1)!\,(a_2-b_2)!}{(b_3-a_3-1)!\,(b_4-a_4-1)!}
\\ &\qquad\times
\iiint\limits_{[0,1]^3}
\frac{\begin{aligned}
x^{a_3-b_1}(1-x)^{b_3-a_3-1}y^{a_4-b_1}(1-y)^{b_4-a_4-1}
\quad\; \\[-3.5pt] \times
z^{a_2-b_1}(1-z)^{a_1-b_2}
\end{aligned}}
{(1-(1-xy)z)^{a_1-b_1+1}}\,\d x\,\d y\,\d z,
\end{align*}
which yields the desired presentation~\eqref{eq:3.7}.
In addition\rom, we mention that the second condition
in~\eqref{eq:3.6} for the parameters~\eqref{eq:3.8}
is equivalent to the condition
\begin{equation}
a_1+a_2+a_3+a_4=b_1+b_2+b_3+b_4-2
\label{eq:3.9}
\end{equation}
for the parameters~\eqref{eq:2.1}.
\end{proof}

The inverse transformation of Rhin--Viola's parameters
to~\eqref{eq:2.1} is defined up to addition of the same
integer to each of the parameters~\eqref{eq:2.1}.
Normalizing the set~\eqref{eq:2.1} by the condition $b_1=1$
we obtain the formulae
\begin{equation}
\begin{alignedat}{4}{4}
a_1&=1+h+q-r, \quad&
a_2&=1+j, \quad&
a_3&=1+h, \quad&
a_4&=1+k,
\\
b_1&=1, \quad&
b_2&=1+h-r, \quad&
b_3&=2+h+l, \quad&
b_4&=2+k+s.
\end{alignedat}
\label{eq:3.10}
\end{equation}
Relations~\eqref{eq:3.8} and~\eqref{eq:3.10} enable us
to describe the action of the generators $\phi,\chi,\vt,\sigma$
of the hypergeometric permutation group~$\varPhi$
from~\cite{RV3} in terms of the parameters~\eqref{eq:2.1}:
\begin{equation}
\begin{aligned}
\phi\colon\biggl(\begin{array}{rrrr}
a_1, & a_2, & a_3, & a_4 \\
  1, & b_2, & b_3, & b_4
\end{array}\biggr)
&\mapsto\biggl(\begin{array}{rrrr}
a_3, & a_2, & a_1, & a_4 \\
  1, & b_2, & b_3, & b_4
\end{array}\biggr),
\\
\chi\colon\biggl(\begin{array}{rrrr}
a_1, & a_2, & a_3, & a_4 \\
  1, & b_2, & b_3, & b_4
\end{array}\biggr)
&\mapsto\biggl(\begin{array}{rrrr}
a_2, & a_1, & a_3, & a_4 \\
  1, & b_2, & b_3, & b_4
\end{array}\biggr),
\\
\vt\colon\biggl(\begin{array}{rrrr}
a_1, & a_2, & a_3, & a_4 \\
  1, & b_2, & b_3, & b_4
\end{array}\biggr)
&\mapsto\biggl(\begin{array}{rr}
b_3-a_1, \qquad a_4, \qquad & \qquad a_2, \;\qquad b_3-a_3 \\
      1, \, b_2+b_3-a_1-a_3, & \, b_3+b_4-a_1-a_3, \, b_3
\end{array}\biggr),
\\
\sigma\colon\biggl(\begin{array}{rrrr}
a_1, & a_2, & a_3, & a_4 \\
  1, & b_2, & b_3, & b_4
\end{array}\biggr)
&\mapsto\biggl(\begin{array}{rrrr}
a_1, & a_2, & a_4, & a_3 \\
  1, & b_2, & b_4, & b_3
\end{array}\biggr).
\end{aligned}
\label{eq:3.11}
\end{equation}
Thus, $\phi,\chi,\sigma$ permute the parameters $a_1,a_2,a_3,a_4$
and $b_3,b_4$ (hence they do not change the quantity~\eqref{eq:3.2}\,),
while the action of the permutation~$\vt$ on the parameters~\eqref{eq:2.1}
is `non-trivial'. In the next section we deduce the group structure
of Rhin and Viola using a classical identity that expresses
Meijer's $G_{4,4}^{2,4}$-function in terms of a well-poised
hypergeometric ${}_7\!F_6$-function. This identity allows us to do
without the integral transform corresponding to~$\vt$ and
to produce another set of generators and another realization
of the same hypergeometric group.

\section{Bailey's identity and the group structure for~$\zeta(3)$}
\label{sec:4}
\begin{proposition}[\textrm{\rm Bailey's identity
\cite{Ba1}, formula~(3.4),
and \cite{Sl}, formula (4.7.1.3)}]
\label{prop:2}
There holds the identity
\begin{equation}
\begin{split}
&
{}_7\!F_6\biggl(\begin{array}{rrrrrrr}
a, & 1+\frac12a, &     b, &     c, &     d, &     e, &     f \\[1pt]
   &   \frac12a, & 1+a-b, & 1+a-c, & 1+a-d, & 1+a-e, & 1+a-f
\end{array}\biggm|1\biggr)
\displaybreak[0]\\ &\qquad
=\frac{\Gamma(1+a-b)\,\Gamma(1+a-c)\,\Gamma(1+a-d)\,
\Gamma(1+a-e)\,\Gamma(1+a-f)}
{\begin{aligned}
\Gamma(1+a)\,\Gamma(b)\,\Gamma(c)\,\Gamma(d)\,
\Gamma(1+a-b-c)\,\Gamma(1+a-b-d)\,
\qquad\quad \\[-3pt] \times
\Gamma(1+a-c-d)\,\Gamma(1+a-e-f)
\end{aligned}}
\\ &\qquad\quad\times
G_{4,4}^{2,4}\biggl(1\biggm|\begin{array}{rrrr}
e+f-a, &       1-b, & 1-c, & 1-d \\
    0, & 1+a-b-c-d, & e-a, & f-a
\end{array}\biggr),
\end{split}
\label{eq:4.1}
\end{equation}
provided that the series on the left-hand side converges.
\end{proposition}

We now set
\begin{equation}
\begin{split}
\wt F(\bh)
&=\wt F(h_0;h_1,h_2,h_3,h_4,h_5)
:=\frac{\Gamma(1+h_0)\cdot\prod_{j=1}^5\Gamma(h_j)}
{\topsmash{\prod_{j=1}^5\Gamma(1+h_0-h_j)}}
\\ &\quad\times
{}_7\!F_6\biggl(\begin{array}{rrrrrr}
h_0, & 1+\frac12h_0, &       h_1, &       h_2, & \dots, &       h_5 \\[1pt]
     &   \frac12h_0, & 1+h_0-h_1, & 1+h_0-h_2, & \dots, & 1+h_0-h_5
\end{array}\biggm|1\biggr)
\end{split}
\label{eq:4.2}
\end{equation}
for the normalized well-poised hypergeometric ${}_7\!F_6$-series.

In the case of integral parameters~$\bh$ satisfying $1+h_0>2h_j$
for each $j=1,\dots,5$, it can be shown that
$\wt F(\bh)$~is a linear form in~$1$ and~$\zeta(3)$
(see, e.g., Section~\ref{sec:8} for the general situation).
Ball's sequence of rational approximations to~$\zeta(3)$
mentioned in Introduction corresponds
to the choice $h_0=3n+2$, $h_1=h_2=h_3=h_4=h_5=n+1$:
\begin{equation}
\begin{gathered}
A_n'\zeta(3)+B_n'
=2n!^2\sum_{t=1}^\infty\biggl(t+\frac n2\biggr)
\frac{(t-1)\dotsb(t-n)\cdot(t+n+1)\dotsb(t+2n)}
{t^4(t+1)^4\dotsb(t+n)^4},
\\
n=1,2,\dots
\end{gathered}
\label{eq:4.3}
\end{equation}
(see~\cite{Ri3}, Section~1.2). Using arguments of Section~\ref{sec:2}
(see also Section~\ref{sec:7} below)
one can show that $D_n\cdot A_n'\in\mathbb Z$
and $D_n^4\cdot B_n'\in\mathbb Z$, which is far from proving
the irrationality of~$\zeta(3)$ since
multiplication of~\eqref{eq:4.3} by~$D_n^4$ leads us to linear
forms with integral coefficients that do not tend to~$0$ as $n\to\infty$.
Rivoal~\cite{Ri3}, Section~5.1, has discovered the coincidence
of Ball's~\eqref{eq:4.3} and Ap\'ery's~\eqref{eq:2.17} sequences
with the use of Zeilberger's {\tt Ekhad} program; the same
result immediately follows from Bailey's identity. Therefore,
one can multiply~\eqref{eq:4.3} by~$D_n^3$ only to obtain
linear forms with integral coefficients! The advantage of
the presentation~\eqref{eq:4.3} of the original Ap\'ery's sequence
consists in the possibility of an `elementary' evaluation of
the series on the right-hand side of~\eqref{eq:4.3} as $n\to\infty$
(see~\cite{Ri3}, Section~5.1, and~\cite{BR} for details).

\begin{lemma}
\label{lem:7}
If condition~\eqref{eq:3.9} holds, then
\begin{align}
&
\frac{\wt G(\ba,\bb)}
{\prod_{j=1}^4(a_j-b_1)!\cdot\prod_{j=1}^4(a_j-b_2)!}
\nonumber\\ &\qquad
=\frac{\wt F(\bh)}
{\prod_{j=1}^5(h_j-1)!\cdot(1+2h_0-h_1-h_2-h_3-h_4-h_5)!},
\label{eq:4.4}
\end{align}
where
\begin{equation}
\begin{gathered}
h_0=b_3+b_4-b_1-a_1=2-2b_1-b_2+a_2+a_3+a_4,
\\
h_1=1-b_1+a_2, \quad h_2=1-b_1+a_3, \quad h_3=1-b_1+a_4,
\\
h_4=b_4-a_1, \quad h_5=b_3-a_1.
\end{gathered}
\label{eq:4.5}
\end{equation}
\end{lemma}

\begin{proof}
Making as before the change of variables $t\mapsto t-b_1+1$ in the contour
integral~\eqref{eq:3.3}, by Lemma~\ref{lem:5} we obtain
$$
\wt G(\ba,\bb)
=G_{4,4}^{2,4}\biggl(1\biggm|\begin{array}{rrrr}
b_1-a_1, & b_1-a_2, & b_1-a_3, & b_1-a_4 \\
      0, & b_1-b_2, & b_1-b_3, & b_1-b_4
\end{array}\biggr).
$$
Therefore, the choice of parameters $h_0,h_1,h_2,h_3,h_4,h_5$ in accordance
with \eqref{eq:4.5} enables us to write down the identity from
Proposition~\ref{prop:2} in the required form~\eqref{eq:4.4}.
\end{proof}

The inverse transformation of the hypergeometric parameters
to~\eqref{eq:2.1} requires a normalization of the parameters~\eqref{eq:2.1}
as in Rhin--Viola's case. Setting $b_1=1$ we obtain
\begin{equation}
\begin{gathered}
a_1=1+h_0-h_4-h_5, \quad
a_2=h_1, \quad a_3=h_2, \quad a_4=h_3,
\\
b_1=1, \quad
b_2=h_1+h_2+h_3-h_0, \quad b_3=1+h_0-h_4, \quad b_4=1+h_0-h_5.
\end{gathered}
\label{eq:4.6}
\end{equation}

We now mention that the permutations~$\fa_{jk}$ of the parameters $a_j,a_k$,
$1\le j<k\le4$, as well as the permutations~$\fb_{12},\fb_{34}$
of the parameters $b_1,b_2$ and $b_3,b_4$ respectively
do not change the quantity on the left-hand side of~\eqref{eq:4.4}.
In a similar way, the permutations~$\fh_{jk}$ of the parameters
$h_j,h_k$, $1\le j<k\le5$, do not change the quantity
on the right-hand side of~\eqref{eq:4.4}. On the other hand,
the permutations~$\fa_{1k}$, $k=2,3,4$, affect nontrivial
transformations of the parameters~$\bh$ and the permutations~$\fh_{jk}$
with $j=1,2,3$ and $k=4,5$ affect nontrivial
transformations of the parameters~$\ba,\bb$. Our nearest goal
is to describe the group~$\fG$ of transformations
of the parameters~\eqref{eq:2.1} and~\eqref{eq:4.5} that is generated
by all (second order) permutations cited above.

\begin{lemma}
\label{lem:8}
The group~$\fG$ can be identified with a subgroup of order~$1920$
of the group~$\fA_{16}$ of even permutations of a $16$-element set\rom;
namely\rom, the group~$\fG$ permutes the parameters
\begin{equation}
c_{jk}=\begin{cases}
a_j-b_k &\mbox{if $a_j\ge b_k$}, \\
b_k-a_j-1 &\mbox{if $a_j<b_k$},
\end{cases}
\qquad j,k=1,2,3,4,
\label{eq:4.7}
\end{equation}
and is generated by following permutations\rom:
\begin{itemize}
\item[\rm(a)] the permutations $\fa_j:=\fa_{j4}$, $j=1,2,3$,
of the $j$th and the fourth lines of the $(4\times4)$-matrix
\begin{equation}
\bc=\pmatrix
c_{11} & c_{12} & c_{13} & c_{14} \\
c_{21} & c_{22} & c_{23} & c_{24} \\
c_{31} & c_{32} & c_{33} & c_{34} \\
c_{41} & c_{42} & c_{43} & c_{44}
\endpmatrix;
\label{eq:4.8}
\end{equation}
\item[\rm(b)] the permutation $\fb:=\fb_{34}$ of the third and the fourth
columns of the matrix~\eqref{eq:4.8}\rom;
\item[\rm(c)] the permutation $\fh:=\fh_{35}$ that has the expression
\begin{equation}
\fh=(c_{11} \; c_{33})(c_{13} \; c_{31})
(c_{22} \; c_{44})(c_{24} \; c_{42})
\label{eq:4.9}
\end{equation}
in terms of the parameters~$\bc$.
\end{itemize}
All these generators have order~$2$.
\end{lemma}

\begin{proof}
The fact that the permutation~$\fh=\fh_{35}$ acts on the
parameters~\eqref{eq:4.7} in accordance with~\eqref{eq:4.9} can be easily
verified with the help of formulae~\eqref{eq:4.5} and~\eqref{eq:4.6}:
\begin{equation}
\fh\colon\biggl(\begin{array}{rrrr}
a_1, & a_2, & a_3, & a_4 \\
  1, & b_2, & b_3, & b_4
\end{array}\biggr)
\mapsto\biggl(\begin{array}{rr}
b_3-a_3, \qquad a_2, \qquad & b_3-a_1, \;\qquad a_4 \qquad \\
      1, \, b_2+b_3-a_1-a_3, & \, b_3, \, b_3+b_4-a_1-a_3
\end{array}\biggr).
\label{eq:4.10}
\end{equation}
As said before, the permutations
$\fa_{jk}$, $1\le j<k\le4$,
and $\fh_{jk}$, $1\le j<k\le5$, belong to the group
$\<\fa_1,\fa_2,\fa_3,\fb,\fh\>$; in addition,
$$
\fb_{12}
=\fh\,\fa_1\,\fa_2\,\fa_1\,\fa_3\,\fh\,
\fb\,\fh\,\fa_3\,\fa_1\,\fa_2\,\fa_1\,\fh.
$$
Therefore, the group~$\fG$ is generated by the elements
in the list~(a)--(c). Obviuosly, these generators
have order~$2$ and belong to~$\fA_{16}$.

We have used a {\tt C++} computer program to find all elements
of the group
\begin{equation}
\fG=\<\fa_1,\fa_2,\fa_3,\fb,\fh\>.
\label{eq:4.11}
\end{equation}
These calculations show that $\fG$~contains exactly
$1920$~permutations. This completes the proof of the lemma.
\end{proof}

\begin{rem}
By Lemma~\ref{lem:8} and relations~\eqref{eq:4.10}
it can be easily verified that
the quantity $b_3+b_4-b_1-b_2$
is stable under the action of~$\fG$.
\end{rem}

Further, a set of parameters~$\bc$, collected in $(4\times4)$-matrix,
is said to be {\it admissible\/} if there exist parameters
$(\ba,\bb)$ such that the elements of the matrix~$\bc$ can be obtained
from them in accordance with~\eqref{eq:4.7} and, moreover,
\begin{equation}
c_{jk}>0 \qquad\mbox{for all}\quad j,k=1,2,3,4.
\label{eq:4.12}
\end{equation}

Comparing the action~\eqref{eq:3.11} of the generators
of the hypergeometric group from \cite{RV3}
on the parameters~\eqref{eq:2.1} with the action of the generators
of the group~\eqref{eq:4.11}, it is easy to see that these two groups
are isomorphic; by~\eqref{eq:4.10} the action of~$\vt$
on~\eqref{eq:2.1} coincides
up to permutations $\fa_1,\fa_2,\fa_3,\fb$ with the action
of~$\fh$. The set of parameters~\eqref{eq:4.7}
is exactly the set~(5.1),~(4.7) from~\cite{RV3}, and
$$
\begin{alignedat}{4}{4}
h&=c_{31}, \quad&
j&=c_{21}, \quad&
k&=c_{41}, \quad&
l&=c_{33},
\\
m&=c_{42}, \quad&
q&=c_{12}, \quad&
r&=c_{32}, \quad&
s&=c_{44}
\end{alignedat}
$$
by~\eqref{eq:3.8}.

On the other hand the hypergeometric group of Rhin and Viola
is embedded into the group~$\fA_{10}$ of even permutations
of a $10$-element set. We can explain this (not so natural, from our
point of view) embedding by pointing out that the following
$10$-element set is stable under~$\fG$:
\begin{alignat*}{2}
h_0-h_1&=b_3+b_4-1-a_1-a_2,
\qquad&
g+h_1&=b_3+b_4-1-a_3-a_4,
\\
h_0-h_2&=b_3+b_4-1-a_1-a_3,
\qquad&
g+h_2&=b_3+b_4-1-a_2-a_4,
\\
h_0-h_3&=b_3+b_4-1-a_1-a_4,
\qquad&
g+h_3&=b_3+b_4-1-a_2-a_3,
\\
h_0-h_4&=b_3-b_1,
\qquad&
g+h_4&=b_4-b_2,
\\
h_0-h_5&=b_4-b_1,
\qquad&
g+h_5&=b_3-b_2,
\end{alignat*}
where $g=1+2h_0-h_1-h_2-h_3-h_4-h_5$.
The matrix~$\bc$ in~\eqref{eq:4.8}
in terms of the parameters~$\bh$ is expressed as
$$
\pmatrix
h_0-h_4-h_5 & g & h_5-1 & h_4-1 \\
h_1-1 & h_0-h_2-h_3 & h_0-h_1-h_4 & h_0-h_1-h_5 \\
h_2-1 & h_0-h_1-h_3 & h_0-h_2-h_4 & h_0-h_2-h_5 \\
h_3-1 & \ h_0-h_1-h_2 \ & \ h_0-h_3-h_4 \ & \ h_0-h_3-h_5 \
\endpmatrix.
$$

The only generator of~$\fG$ in the list~(a)--(c)
that acts nontrivially on the parameters~$\bh$ is
the permutation~$\fa_1$. Its action is
\begin{align*}
&
(h_0;h_1,h_2,h_3,h_4,h_5)
\mapsto(1+2h_0-h_3-h_4-h_5;
\\ &\qquad
h_1,h_2,1+h_0-h_4-h_5,1+h_0-h_3-h_5,1+h_0-h_3-h_4),
\end{align*}
and we have discovered the corresponding hypergeometric
${}_7\!F_6$-identity in \cite{Ba2}, formula~(2.2).

The subgroup $\fG_1$ of~$\fG$ generated by the permutations
$\fa_{jk}$, $1\le j<k\le4$, and $\fb_{12},\fb_{34}$,
has order $4!\cdot2!\cdot2!=96$. The quantity $\wt G(\ba,\bb)$
is stable under the action of this group, hence we can present
the group action on the parameters by indicating
$1920/96=20$ representatives of left cosets
$\fG/\fG_1=\{\fq_j\fG_1,\ j=1,\dots,20\}$;
namely,
\begin{alignat*}{4}
\fq_1&=\id,
\;\;&
\fq_2&=\fa_1\,\fa_2\,\fa_3\,\fh,
\;\;&
\fq_3&=\fa_1\,\fh,
\;\;&
\fq_4&=\fa_2\,\fa_1\,\fh,
\\
\fq_5&=\fh,
\;\;&
\fq_6&=\fh\,\fa_1\,\fa_2\,\fa_3\,\fh,
\;\;&
\fq_7&=\fa_2\,\fa_3\,\fh,
\;\;&
\fq_8&=\fa_3\,\fh,
\\
\fq_9&=\fh\,\fa_3\,\fb\,\fh,
\;\;&
\fq_{10}&=\fa_1\,\fa_2\,\fh\,\fa_1\,\fa_2\,\fb\,\fh,
\;\;&
\fq_{11}&=\fa_2\,\fh\,\fa_3\,\fa_2\,\fb\,\fh,
\;\;&
\fq_{12}&=\fb\,\fh,
\\
\fq_{13}&=\fa_2\,\fa_3\,\fb\,\fh,
\;\;&
\fq_{14}&=\fa_3\,\fb\,\fh,
\;\;&
\fq_{15}&=\fa_1\,\fa_2\,\fa_3\,\fb\,\fh,
\;\;&
\fq_{16}&=\fa_1\,\fb\,\fh,
\\
\fq_{17}&=\fa_2\,\fa_1\,\fb\,\fh,
\;\;&
\fq_{18}&=\fa_2\,\fh\,\fa_1\,\fa_2\,\fb\,\fh,
\;\;&
\fq_{19}&=\fa_3\,\fh\,\fa_1\,\fb\,\fh,
\;\;&
\fq_{20}&=\fh\,\fa_1\,\fb\,\fh;
\end{alignat*}
we choose the representatives with the shortest presentation
in terms of the generators from the list~(a)--(c).
The images of any set of parameters $(\ba,\bb)$
under the action of these representatives
can be normalized by the condition
$b_1=1$ and ordered in accordance with~\eqref{eq:2.9}.
We also point out that the group~$\fG_1$ contains the subgroup
$\fG_0=\<\fa_{12}\fb_{12},\fa_{34}\fb_{34}\>$
of order~$4$, which does not change the quantity
$G(\ba,\bb)$. This fact shows us that for fixed data
$(\ba,\bb)$ only the $480$~elements
$\fq_j\fa$, where $j=1,\dots,20$ and $\fa\in\fS_4$~is an arbitrary
permutation of the parameters $a_1,a_2,a_3,a_4$,
produce `perceptable' actions on the quantity~\eqref{eq:2.7}.
Hence we will restrict ourselves to the consideration of only these
$480$~permutations from~$\fG/\fG_0$.

In the same way one can consider the subgroup
$\fG_1'\subset\fG$ of order $5!=120$ generated by the permutations
$\fh_{jk}$, $1\le j<k\le5$. This group acts trivially
on the quantity $\wt F(\bh)$. The corresponding
$1920/120=16$ representatives of left cosets $\fG/\fG_1'$
can be chosen so that for the images of the set of parameters~$\bh$
we have
$$
1\le h_1\le h_2\le h_3\le h_4\le h_5;
$$
of course $h_0>2h_5$.

For an admissible set of parameters~\eqref{eq:4.7} consider
the quantity
\begin{equation}
H(\bc):=G(\ba,\bb)
=\frac{c_{33}!\,c_{44}!}{c_{11}!\,c_{22}!}\wt G(\ba,\bb).
\label{eq:4.13}
\end{equation}
Since the group~$\fG$ does not change~\eqref{eq:4.4},
we arrive at the following statement.

\begin{lemma}[\textrm{\rm cf.~\cite{RV3}, Section~4}]
\label{lem:9}
The quantity
\begin{equation}
\frac{H(\bc)}{\varPi(\bc)},
\qquad\mbox{where}\quad
\varPi(\bc)=c_{21}!\,c_{31}!\,c_{41}!\,
c_{12}!\,c_{32}!\,c_{42}!\,c_{33}!\,c_{44}!\,,
\label{eq:4.14}
\end{equation}
is stable under the action of~$\fG$.
\end{lemma}

\section{Irrationality measure of Rhin and Viola for~$\zeta(3)$}
\label{sec:5}
Throught this section the set of parameters~\eqref{eq:2.1}
will depend on a positive integer~$n$ in the following way:
\begin{equation}
\begin{alignedat}{4}{4}
a_1&=\alpha_1n+1, \quad&
a_2&=\alpha_2n+1, \quad&
a_3&=\alpha_3n+1, \quad&
a_4&=\alpha_4n+1,
\\
b_1&=\beta_1n+1, \quad&
b_2&=\beta_2n+1, \quad&
b_3&=\beta_3n+2, \quad&
b_4&=\beta_4n+2,
\end{alignedat}
\label{eq:5.1}
\end{equation}
where the {\it new\/} integral parameters
(`directions') $(\balpha,\bbeta)$
satisfy by~\eqref{eq:2.2}, \eqref{eq:3.9}, and~\eqref{eq:4.12}
the following conditions:
\begin{gather}
\{\beta_1,\beta_2\}
<\{\alpha_1,\alpha_2,\alpha_3,\alpha_4\}
<\{\beta_3,\beta_4\},
\label{eq:5.2}
\\
\alpha_1+\alpha_2+\alpha_3+\alpha_4
=\beta_1+\beta_2+\beta_3+\beta_4.
\label{eq:5.3}
\end{gather}
The version of the set $(\balpha,\bbeta)$ ordered as in~\eqref{eq:2.9}
is denoted by $(\balpha^*,\bbeta^*)$.

To the parameters $(\balpha,\bbeta)$ we assign the admissible
$(4\times4)$-matrix~$\bc$ with elements
\begin{equation}
c_{jk}=\begin{cases}
\alpha_j-\beta_k & \mbox{if $\alpha_j>\beta_k$}, \\
\beta_k-\alpha_j & \mbox{if $\alpha_j<\beta_k$},
\end{cases}
\qquad j,k=1,2,3,4,
\label{eq:5.4}
\end{equation}
hence the set of parameters $\bc\cdot n$ corresponds to~\eqref{eq:5.1}.
With any admissible matrix~$\bc$ we relate the following characteristics:
$$
\begin{gathered}
\begin{aligned}
m_0=m_0(\bc)
&:=\max_{1\le j,k\le4}\{c_{jk}\}>0,
\\
m_1=m_1(\bc)
&:=\beta_4^*-\alpha_1^*=\max_{1\le j\le4}\{c_{j3},c_{j4}\},
\\
m_2=m_2(\bc)
&:=\max\{\alpha_1-\beta_1,\alpha_2-\beta_2,
\beta_4^*-\alpha_3,\beta_4^*-\alpha_4,\beta_3^*-\alpha_1^*\}
\\ &\phantom:
=\max\{c_{11},c_{1k},c_{22},c_{2k},c_{34},c_{44},c_{33},c_{43}\},
\end{aligned}
\\
\mbox{where}\quad
k=\begin{cases}
3 &\mbox{if $\beta_4=\beta_4^*$ (i.e., $c_{13}\le c_{14}$)}, \\
4 &\mbox{if $\beta_3=\beta_4^*$ (i.e., $c_{13}\ge c_{14}$)},
\end{cases}
\end{gathered}
$$
and write the claim of Lemma~\ref{lem:4}
for the quantity~\eqref{eq:4.13} as
\begin{equation}
D_{m_1(\bc)n}^2\cdot D_{m_2(\bc)n}\cdot H(\bc n)
\in2\mathbb Z\zeta(3)+\mathbb Z.
\label{eq:5.5}
\end{equation}

Fix now a set of directions $(\balpha,\bbeta)$
satisfying conditions \eqref{eq:5.2}, \eqref{eq:5.3}, and the
corresponding set of parameters~\eqref{eq:5.4}. In view of the results
of Section~\ref{sec:4}, we will consider the set
$\sM_0=\sM_0(\balpha,\bbeta)=\sM_0(\bc)$
of $20$~ordered collections $(\balpha',\bbeta')$
corresponding to $\fq_j(\balpha,\bbeta)$, $j=1,\dots,20$,
and the set $\sM=\sM(\balpha,\bbeta)=\sM(\bc):=\{\fa\sM_0\}$
of $480$~such collections, where $\fa\in\fS_4$~is an arbitrary
permutation of the parameters
$\alpha_1,\alpha_2,\alpha_3,\alpha_4$
(equivalently, of the lines of the matrix~$\bc$).
To each prime number~$p$ we assign the exponent
$$
\nu_p=\max_{\bc'\in\sM}
\ord_p\frac{\varPi(\bc n)}{\varPi(\bc'n)}
$$
and consider the quantity
\begin{equation}
\Phi_n=\Phi_n(\bc):=\prod_{\sqrt{m_0n}<p\le m_3n}p^{\nu_p},
\label{eq:5.6}
\end{equation}
where $m_3=m_3(\bc):=\min\{m_1(\bc),m_2(\bc)\}$.

\begin{lemma}
\label{lem:10}
For any positive integer~$n$ there holds the inclusion
$$
D_{m_1n}^2\cdot D_{m_2n}
\cdot\Phi_n^{-1}\cdot H(\bc n)
\in2\mathbb Z\zeta(3)+\mathbb Z.
$$
\end{lemma}

\begin{proof}
The inclusions
\begin{equation}
D_{m_1n}^2\cdot D_{m_2n}
\cdot\Phi_n^{-1}\cdot H(\bc n)
\in2\mathbb Z_p\zeta(3)+\mathbb Z_p
\label{eq:5.7}
\end{equation}
for $p\le\sqrt{m_0n}$ and $p>m_3n$ follow from~\eqref{eq:5.5}
since $\ord_p\Phi_n^{-1}=0$.

Using the stability of the quantity~\eqref{eq:4.14}
under the action of any permutation from the group~$\fG$,
by~\eqref{eq:5.5} we deduce that
\begin{align*}
&
D_{m_1(\bc')n}^2\cdot D_{m_2(\bc')n}
\cdot\frac{\varPi(\bc'n)}{\varPi(\bc n)}\cdot H(\bc n)
\\ &\qquad
=D_{m_1(\bc')n}^2\cdot D_{m_2(\bc')n}\cdot H(\bc'n)
\in2\mathbb Z\zeta(3)+\mathbb Z,
\qquad \bc'\in\sM,
\end{align*}
which yields the inclusions~\eqref{eq:5.7} for the primes~$p$
in the interval $\sqrt{m_0n}<p\le m_3n$ since
\begin{align*}
\ord_p\bigl(D_{m_1(\bc')n}^2\cdot D_{m_2(\bc')n}\bigr)
\le3
&=\ord_p\bigl(D_{m_3(\bc)n}^3\bigr)
\\
&=\ord_p\bigl(D_{m_1(\bc)n}^2\cdot D_{m_2(\bc)n}\bigr),
\qquad \bc'\in\sM(\bc)
\end{align*}
in this case. The proof is complete.
\end{proof}

The asymptotics of the numbers $D_{m_1n},D_{m_2n}$ in~\eqref{eq:5.7}
is determined from the prime number theorem:
$$
\lim_{n\to\infty}\frac{\log D_{m_jn}}n=m_j,
\qquad j=1,2.
$$
For the study of the asymptotic behaviour of~\eqref{eq:5.6}
as $n\to\infty$ we introduce the function
\begin{align*}
\phi(x)=\smash{\max_{\bc'\in\sM}}&\bigl(
\[c_{21}x\]+\[c_{31}x\]+\[c_{41}x\]+\[c_{12}x\]
\\ &\qquad
+\[c_{32}x\]+\[c_{42}x\]+\[c_{33}x\]+\[c_{44}x\]
\\ &\qquad
-\[c_{21}'x\]-\[c_{31}'x\]-\[c_{41}'x\]-\[c_{12}'x\]
\\ &\qquad
-\[c_{32}'x\]-\[c_{42}'x\]-\[c_{33}'x\]-\[c_{44}'x\]\bigr),
\end{align*}
where $\[\,\cdot\,\]$~is the integral part of a number.
Then $\nu_p=\phi(n/p)$ since
$\ord_pN!=\[N/p\]$ for any
integer~$N$ and any prime $p>\sqrt N$.

Note that the function $\phi(x)$ is periodic (with period~$1$)
since
\begin{align*}
&
c_{21}+c_{31}+c_{41}+c_{12}+c_{32}+c_{42}+c_{33}+c_{44}
=2(\beta_3+\beta_4-\beta_1-\beta_2)
\\ &\qquad
=c_{21}'+c_{31}'+c_{41}'+c_{12}'+c_{32}'+c_{42}'+c_{33}'+c_{44}'
\end{align*}
(see Remark to Lemma~\ref{lem:8});
moreover, the function $\phi(x)$
takes only non-negative integral values.

\begin{lemma}
\label{lem:11}
The number~\eqref{eq:5.6} satisfies the limit relation
\begin{equation}
\lim_{n\to\infty}\frac{\log\Phi_n}n
=\int_0^1\phi(x)\,\d\psi(x)-\int_0^{1/m_3}\phi(x)\frac{\d x}{x^2},
\label{eq:5.8}
\end{equation}
where $\psi(x)$~is the logarithmic derivative of the gamma function.
\end{lemma}

\begin{proof}
This result follows from the arithmetic scheme
of Chudnovsky--Rukha\-dze--Hata
and is based on the above-cited properties of the function~$\phi(x)$
(see \cite{Zu3}, Lemma~4.4). Subtraction on the right-hand side
of~\eqref{eq:5.8} `removes' the primes $p>m_3n$ that do not enter
the product~$\Phi_n$ in~\eqref{eq:5.6}.
\end{proof}

The asymptotic behaviour of linear forms
$$
H_n:=H(\bc n)=2A_n\zeta(3)-B_n
$$
and their coefficients $A_n,B_n$ can be deduced from
Lemma~\ref{lem:6} and~\cite{RV3}, the arguments before Theorem~5.1;
another `elementary' way is based on the presentation
\begin{align}
H(\bc)
&=\frac{(h_0-h_1-h_2)!\,(h_0-h_1-h_3)!\,(h_0-h_2-h_4)!\,(h_0-h_3-h_5)!}
{(h_4-1)!\,(h_5-1)!}
\nonumber\\ &\qquad\times 
\wt F(\bh)
\label{eq:5.9}
\end{align}
and the arguments of Ball (see \cite{BR} or \cite{Ri3}, Section~5.1).
But the same asymptotic problem can be solved directly
on the basis of Lemma~\ref{lem:5} with the use of the asymptotics of the gamma
function and the saddle-point method. We refer the reader to~\cite{Ne1}
and~\cite{Zu3}, Sections~2 and~3, for details of this approach;
here we only state the final result.

\begin{lemma}
\label{lem:12}
Let $\tau_0<\tau_1$ be the \rom(real\rom) zeros of the quadratic polynomial
$$
(\tau-\alpha_1)(\tau-\alpha_2)(\tau-\alpha_3)(\tau-\alpha_4)
-(\tau-\beta_1)(\tau-\beta_2)(\tau-\beta_3)(\tau-\beta_4)
$$
\rom(it can be easily verified that $\beta_2^*<\tau_0<\alpha_1^*$
and $\tau_1>\alpha_4^*$\rom)\rom;
the function $f_0(\tau)$ in the cut $\tau$-plane
$\mathbb C\setminus(-\infty,\beta_2^*]\cup[\alpha_1^*,+\infty)$
is given by the formula
\begin{align*}
f_0(\tau)
&=\alpha_1\log(\alpha_1-\tau)
+\alpha_2\log(\alpha_2-\tau)
+\alpha_3\log(\alpha_3-\tau)
+\alpha_4\log(\alpha_4-\tau)
\\ &\quad
-\beta_1\log(\tau-\beta_1)
-\beta_2\log(\tau-\beta_2)
-\beta_3\log(\beta_3-\tau)
-\beta_4\log(\beta_4-\tau)
\\ &\quad
-(\alpha_1-\beta_1)\log(\alpha_1-\beta_1)
-(\alpha_2-\beta_2)\log(\alpha_2-\beta_2)
\\ &\quad
+(\beta_3-\alpha_3)\log(\beta_3-\alpha_3)
+(\beta_4-\alpha_4)\log(\beta_4-\alpha_4),
\end{align*}
where the logarithms
take real values for real $\tau\in(\beta_2^*,\alpha_1^*)$.
Then
$$
\lim_{n\to\infty}\frac{\log|H_n|}n=f_0(\tau_0),
\qquad
\limsup_{n\to\infty}\frac{\log\max\{|A_n|,|B_n|\}}n
\le\Re f_0(\tau_1).
$$
\end{lemma}

Combining results of Lemmas~\ref{lem:11} and~\ref{lem:12},
as in~\cite{RV3}, Theorem~5.1,
we deduce the following statement.

\begin{proposition}
\label{prop:3}
In the above notation let
$$
\begin{gathered}
C_0=-f_0(\tau_0), \qquad C_1=\Re f_0(\tau_1),
\\
C_2=2m_1+m_2
-\biggl(\int_0^1\phi(x)\,\d\psi(x)
-\int_0^{1/m_3}\phi(x)\frac{\d x}{x^2}\biggr).
\end{gathered}
$$
If $C_0>C_2$\rom, then
$$
\mu(\zeta(3))\le\frac{C_0+C_1}{C_0-C_2}.
$$
\end{proposition}

Looking over all integral directions $(\balpha,\bbeta)$
satisfying the relation
\begin{equation}
\alpha_1+\alpha_2+\alpha_3+\alpha_4
=\beta_1+\beta_2+\beta_3+\beta_4
\le200
\label{eq:5.10}
\end{equation}
by means of a program for the calculator {\tt GP-PARI} we have
discovered that the best estimate for~$\mu(\zeta(3))$
is given by Rhin and Viola in~\cite{RV3}.

\begin{theorem}[\textrm{\rm\cite{RV3}}]
\label{th:1}
The irrationality exponent of~$\zeta(3)$ satisfies the estimate
\begin{equation}
\mu(\zeta(3))\le5.51389062\dots\,.
\label{eq:5.11}
\end{equation}
\end{theorem}

\begin{proof}
The optimal set of directions $(\balpha,\bbeta)$
(up to the action of~$\fG$) is as follows:
\begin{equation}
\begin{alignedat}{4}{4}
\alpha_1&=18, \quad& \alpha_2&=17, \quad&
\alpha_3&=16, \quad& \alpha_4&=19,
\\
\beta_1&=0, \quad& \beta_2&=7, \quad&
\beta_3&=31, \quad& \beta_4&=32.
\end{alignedat}
\label{eq:5.12}
\end{equation}
Then,
\begin{alignat*}{2}
\tau_0&=8.44961969\dots,
\qquad&
C_0&=-f_0(\tau_0)=47.15472079\dots,
\\
\tau_1&=27.38620119\dots,
\qquad&
C_1&=\Re f_0(\tau_0)=48.46940964\dots\,.
\end{alignat*}

The set $\sM_0$ in this case consists of the following elements:
\begin{gather*}
\biggl(\begin{array}{rrrr}
16, & 17, & 18, & 19 \\0, & 7, & 31, & 32
\end{array}\biggr),
\;
\biggl(\begin{array}{rrrr}
12, & 14, & 16, & 18 \\ 0, & 2, & 27, & 31
\end{array}\biggr),
\;
\biggl(\begin{array}{rrrr}
12, & 15, & 17, & 18 \\ 0, & 3, & 28, & 31
\end{array}\biggr),
\;
\biggl(\begin{array}{rrrr}
14, & 15, & 18, & 19 \\ 0, & 5, & 30, & 31
\end{array}\biggr),
\\
\biggl(\begin{array}{rrrr}
13, & 15, & 17, & 19 \\ 0, & 4, & 29, & 31
\end{array}\biggr),
\;
\biggl(\begin{array}{rrrr}
13, & 14, & 15, & 16 \\ 0, & 1, & 25, & 32
\end{array}\biggr),
\;
\biggl(\begin{array}{rrrr}
13, & 14, & 16, & 19 \\ 0, & 3, & 28, & 31
\end{array}\biggr),
\;
\biggl(\begin{array}{rrrr}
12, & 13, & 16, & 17 \\ 0, & 1, & 26, & 31
\end{array}\biggr),
\\
\biggl(\begin{array}{rrrr}
11, & 14, & 15, & 18 \\ 0, & 1, & 27, & 30
\end{array}\biggr),
\;
\biggl(\begin{array}{rrrr}
11, & 15, & 16, & 18 \\ 0, & 2, & 28, & 30
\end{array}\biggr),
\;
\biggl(\begin{array}{rrrr}
12, & 13, & 14, & 19 \\ 0, & 1, & 28, & 29
\end{array}\biggr),
\;
\biggl(\begin{array}{rrrr}
14, & 16, & 17, & 19 \\ 0, & 5, & 29, & 32
\end{array}\biggr),
\\
\biggl(\begin{array}{rrrr}
14, & 15, & 16, & 19 \\ 0, & 4, & 28, & 32
\end{array}\biggr),
\;
\biggl(\begin{array}{rrrr}
13, & 14, & 16, & 17 \\ 0, & 2, & 26, & 32
\end{array}\biggr),
\;
\biggl(\begin{array}{rrrr}
13, & 15, & 16, & 18 \\ 0, & 3, & 27, & 32
\end{array}\biggr),
\;
\biggl(\begin{array}{rrrr}
13, & 16, & 17, & 18 \\ 0, & 4, & 28, & 32
\end{array}\biggr),
\\
\biggl(\begin{array}{rrrr}
15, & 16, & 18, & 19 \\ 0, & 6, & 30, & 32
\end{array}\biggr),
\;
\biggl(\begin{array}{rrrr}
12, & 15, & 16, & 19 \\ 0, & 3, & 29, & 30
\end{array}\biggr),
\;
\biggl(\begin{array}{rrrr}
12, & 14, & 15, & 19 \\ 0, & 2, & 28, & 30
\end{array}\biggr),
\;
\biggl(\begin{array}{rrrr}
10, & 15, & 16, & 17 \\ 0, & 1, & 28, & 29
\end{array}\biggr);
\end{gather*}
an easy verification shows that $m_1=m_3=16$ and $m_2=18$. The function
$\phi(x)$ for $x\in[0,1)$ is defined by the formula
$$
\phi(x)=\begin{cases}
0 & \mbox{if $x\in[0,1)\setminus\varOmega_E$}, \\
1 & \mbox{if $x\in\varOmega_E\setminus\varOmega_E'$}, \\
2 & \mbox{if $x\in\varOmega_E'$},
\end{cases}
$$
where the sets $\varOmega_E$ and $\varOmega_E'$ are indicated
in~\cite{RV3}, p.~292. Hence
\begin{align*}
C_2
&=2m_1+m_2
-\biggl(\int_0^1\phi(x)\,\d\psi(x)
-\int_0^{1/m_3}\phi(x)\frac{\d x}{x^2}\biggr)
\\
&=2\cdot16+18-(24.18768530\hdots-4)
=29.81231469\dots,
\end{align*}
and by Proposition~\ref{prop:3}
we obtain the required estimate~\eqref{eq:5.11}.
\end{proof}

Note that the choice~\eqref{eq:5.12} gives us the function~$\phi(x)$
ranging in the set $\{0,1,2\}$; any other element of~$\sM$
produces the same estimate of the irrationality exponent~\eqref{eq:5.11}
with $\phi(x)$ ranging in $\{0,1,2,3\}$.

The previous record
\begin{equation}
\mu(\zeta(3))\le7.37795637\dots
\label{eq:5.13}
\end{equation}
due to Hata~\cite{Ha5} can be achieved by the choice of the parameters
\begin{equation}
\begin{alignedat}{4}{4}
\alpha_1&=8, \quad& \alpha_2&=7, \quad&
\alpha_3&=8, \quad& \alpha_4&=9,
\\
\beta_1&=0, \quad& \beta_2&=1, \quad&
\beta_3&=15, \quad& \beta_4&=16,
\end{alignedat}
\label{eq:5.14}
\end{equation}
and the action of the group $\fG_1/\fG_0$ of order just~$4!=24$
(we can regard this as a $(\ba,\bb)$-{\it trivial action\/}).
For directions $(\balpha,\bbeta)$ satisfying the relation
$$
\alpha_1+\alpha_2+\alpha_3+\alpha_4
\le\beta_1+\beta_2+\beta_3+\beta_4
\le200
$$
(instead of~\eqref{eq:5.10}\,) we have verified that the choice~\eqref{eq:5.14}
corresponding to Hata's case produces the best estimate of the
irrationality exponent for~$\zeta(3)$ in the class of
$(\ba,\bb)$-trivial actions.
In that case we are able to use the inequality
$$
\alpha_1+\alpha_2+\alpha_3+\alpha_4
\le\beta_1+\beta_2+\beta_3+\beta_4
$$
instead of~\eqref{eq:5.3} since we do not use Bailey's identity.
The mysterious thing is that the action of the full group~$\fG$
does not produce a better result than \eqref{eq:5.13}
for the parameters~\eqref{eq:5.14}.

\section{Overview of the group structure for $\zeta(2)$}
\label{sec:6}
To a set of integral parameters
\begin{equation}
(\ba,\bb)
=\biggl(\begin{array}{rrr}
a_1, & a_2, & a_3 \\
b_1, & b_2, & b_3
\end{array}\biggr)
\label{eq:6.1}
\end{equation}
satisfying the conditions
\begin{gather}
\{b_1\}\le\{a_1,a_2,a_3\}<\{b_2,b_3\},
\nonumber\\
a_1+a_2+a_3\le b_1+b_2+b_3-2,
\label{eq:6.1a}
\end{gather}
we assign the rational function
\begin{align*}
R(t)
=R(\ba,\bb;t)
&:=\frac{(b_2-a_2-1)!\,(b_3-a_3-1)!}{(a_1-b_1)!}
\\ &\phantom:\qquad\times 
\frac{\Gamma(t+a_1)\,\Gamma(t+a_2)\,\Gamma(t+a_3)}
{\Gamma(t+b_1)\,\Gamma(t+b_2)\,\Gamma(t+b_3)}
\\ &\phantom:
=\prod_{j=1}^3R_j(t),
\end{align*}
where the functions $R_1(t),R_2(t)$, and $R_3(t)$ are defined
in~\eqref{eq:2.5}. Condition~\eqref{eq:6.1a} yields \eqref{eq:2.6},
hence the (hypergeometric) series
\begin{equation}
G(\ba,\bb):=\sum_{t=t_0}^\infty R(t)
\qquad\mbox{with}\quad
1-\min\{a_1,a_2,a_3\}\le t_0\le1-b_1
\label{eq:6.2}
\end{equation}
is well-defined. Expanding the rational function $R(t)$
in a sum of partial fractions
and applying Lemmas~\ref{lem:1} and~\ref{lem:3}
we arrive at the following assertion.

\begin{lemma}[\textrm{\rm cf.~Lemma~\ref{lem:4}}]
\label{lem:13}
The quantity~\eqref{eq:6.2} is a rational form in~$1$ and~$\zeta(2)$
with rational coefficients\rom:
\begin{equation}
G(\ba,\bb)=A\zeta(2)-B;
\label{eq:6.3}
\end{equation}
in addition\rom,
$$
A\in\mathbb Z, \qquad
D_{b_3^*-a_1^*-1}
\cdot D_{\max\{a_1-b_1,b_3^*-a_2-1,b_3^*-a_3-1,b_2^*-a_1^*-1\}}
\cdot B\in\mathbb Z,
$$
where $(\ba^*,\bb^*)$ is the ordered version of the set~\eqref{eq:6.1}\rom:
\begin{equation}
\begin{gathered}
\{b_1^*\}=\{b_1\},
\quad
\{a_1^*,a_2^*,a_3^*\}=\{a_1,a_2,a_3\},
\quad
\{b_2^*,b_3^*\}=\{b_2,b_3\},
\\
b_1^*\le a_1^*\le a_2^*\le a_3^*<b_2^*\le b_3^*.
\end{gathered}
\label{eq:6.4}
\end{equation}
\end{lemma}

By Proposition~\ref{prop:1} the series~\eqref{eq:6.2} can be written
as the double real integral
$$
G(\ba,\bb)
=\iint\limits_{[0,1]^2}
\frac{x^{a_2-b_1}(1-x)^{b_2-a_2-1}y^{a_3-b_1}(1-y)^{b_3-a_3-1}}
{(1-xy)^{a_1-b_1+1}}\,\d x\,\d y,
$$
hence we can identify the quantity~\eqref{eq:6.2} with
the corresponding integral $I(h,i,j,k,l)$ from~\cite{RV2}
by setting
$$
\begin{gathered}
h=a_2-b_1, \quad i=b_2-a_2-1, \quad j=b_3-a_3-1,
\\
k=a_3-b_1, \quad l=(b_1+b_2+b_3-2)-(a_1+a_2+a_3);
\end{gathered}
$$
the inverse transformation (after the normalization $b_1=1$) is as
follows:
$$
\begin{alignedat}{3}{3}
a_1&=1+i+j-l, \quad& a_2&=1+h, \quad& a_3&=1+k,
\\
b_1&=1, \quad& b_2&=2+h+i, \quad& b_3&=2+j+k.
\end{alignedat}
$$
In the further discussion we keep the normalization $b_1=1$.

The series
$$
\wt G(\ba,\bb)
:=\frac{\Gamma(a_1)\Gamma(a_2)\Gamma(a_3)}
{\Gamma(b_1)\Gamma(b_2)\Gamma(b_3)}
\cdot{}_3\!F_2\biggl(\begin{array}{rrr}
a_1, & a_2, & a_3 \\
     & b_2, & b_3
\end{array}\biggm|1\biggr)
$$
and
\begin{align*}
\wt F(\bh)
&=\wt F(h_0;h_1,h_2,h_3,h_4)
:=\frac{\Gamma(1+h_0)\cdot\prod_{j=1}^4\Gamma(h_j)}
{\prod_{j=1}^4\Gamma(1+h_0-h_j)}
\\ &\quad\times
{}_6\!F_5\biggl(\begin{array}{rrrrr}
h_0, & 1+\frac12h_0, &       h_1, & \dots, &       h_4 \\[1pt]
     &   \frac12h_0, & 1+h_0-h_1, & \dots, & 1+h_0-h_4
\end{array}\biggm|-1\biggr)
\end{align*}
play the same role as~\eqref{eq:3.2} and~\eqref{eq:4.2} played before
since one has
\begin{align}
&
\frac{\wt G(\ba,\bb)}
{\Gamma(a_1)\,\Gamma(a_2)\,\Gamma(a_3)\,\Gamma((b_2+b_3)-(a_1+a_2+a_3))}
\nonumber\\ &\qquad
=\frac{\wt F(\bh)}
{\Gamma(h_1)\,\Gamma(h_2)\,\Gamma(h_3)\,\Gamma(h_4)}
\label{eq:6.5}
\end{align}
where
$$
\begin{gathered}
h_0=b_2+b_3-1-a_1, \quad h_1=a_2, \quad h_2=a_3,
\\
h_3=b_3-a_1, \quad h_4=b_2-a_1,
\end{gathered}
$$
and
$$
\begin{alignedat}{3}{3}
a_1&=1+h_0-h_3-h_4, \quad& a_2&=h_1, \quad& a_3&=h_2,
\\
b_1&=1, \quad& b_2&=1+h_0-h_3, \quad& b_3&=1+h_0-h_4,
\end{alignedat}
$$
by Whipple's identity~\cite{Ba3}, Section~4.4, formula~(2).
The permutations $\fa_{jk}$, $1\le j<k\le3$, of the parameters $a_j,a_k$,
the permutation~$\fb_{23}$ of~$b_2,b_3$, and the permutations
$\fh_{jk}$, $1\le j<k\le4$, of the parameters $h_j,h_k$
do not change the quantity~\eqref{eq:6.5}. Hence we can consider
the group~$\fG$ generated by these permutations and naturally embed it
into the group~$\fS_{10}$ of permutations of the $10$-element set
\begin{align*}
c_{00}&=(b_2+b_3)-(a_1+a_2+a_3)-1,
\\
c_{jk}&=\begin{cases}
a_j-b_k &\mbox{if $a_j\ge b_k$}, \\
b_k-a_j-1 &\mbox{if $a_j<b_k$},
\end{cases}
\qquad j,k=1,2,3.
\end{align*}
The group~$\fG$ is generated by the permutations $\fa_1:=\fa_{13}$,
$\fa_2:=\fa_{23}$, $\fb:=\fb_{23}$, which can be regarded
as permutations of lines and columns of the `$(4\times4)$-matrix'
\begin{equation}
\bc=\pmatrix
c_{00} &        &        &        \\
       & c_{11} & c_{12} & c_{13} \\
       & c_{21} & c_{22} & c_{23} \\
       & c_{31} & c_{32} & c_{33}
\endpmatrix,
\label{eq:6.6}
\end{equation}
and the $(\ba,\bb)$-nontrivial permutation $\fh:=\fh_{23}$,
$$
\fh=(c_{00} \; c_{22})(c_{11} \; c_{33})(c_{13} \; c_{31});
$$
these four generators have order~$2$. It can be easily verified
that the group $\fG=\<\fa_1,\fa_2,\fb,\fh\>$
has order~$120$;
in fact, we require only the $60$~representatives of~$\fG/\fG_0$,
where the group $\fG_0=\{\id,\fa_{23}\fb_{23}\}$ acts trivially
on the quantity
$$
H(\bc):=G(\ba,\bb)
=\frac{c_{22}!\,c_{33}!}{c_{11}!}\wt G(\ba,\bb).
$$
Thus, we can summarize the above as follows.

\begin{lemma}[\textrm{\rm cf.~\cite{RV2}, Section~3}]
\label{lem:14}
The quantity
$$
\frac{H(\bc)}{\varPi(\bc)},
\qquad\mbox{where}\quad
\varPi(\bc)=c_{00}!\,c_{21}!\,c_{31}!\,c_{22}!\,c_{33}!\,,
$$
is stable under the action of $\fG=\<\fa_1,\fa_2,\fb,\fh\>$.
\end{lemma}

If one shifts indices of $c_{jk}$ by one then the group~$\fG$
for~$\zeta(2)$ can be naturally regarded as a subgroup of
the group~$\fG$ for~$\zeta(3)$ (compare the generators of both groups).
The group~$\fG$ for~$\zeta(2)$ coincides with
the group~$\boldsymbol\Phi$ of Rhin and Viola from~\cite{RV2}
since permutations $\phi,\sigma\in\boldsymbol\Phi$ are
$(\ba,\bb)$-trivial in our terms and for $\tau\in\boldsymbol\Phi$
we have
$$
\tau=\fa_2\,\fa_1\,\fb\,\fh\,\fa_2\,\fa_1\,\fb\,\fh.
$$

We now fix an arbitrary positive integer~$n$ and
integral directions $(\balpha,\bbeta)$ satisfying the conditions
\begin{gather*}
\{\beta_1=0\}
<\{\alpha_1,\alpha_2,\alpha_3\}
<\{\beta_2,\beta_3\},
\\
\alpha_1+\alpha_2+\alpha_3
\le\beta_1+\beta_2+\beta_3,
\end{gather*}
so that the parameters~\eqref{eq:6.1} are expressed as follows:
\begin{equation}
\begin{alignedat}{3}{3}
a_1&=\alpha_1n+1, \quad&
a_2&=\alpha_2n+1, \quad&
a_3&=\alpha_3n+1,
\\
b_1&=\beta_1n+1, \quad&
b_2&=\beta_2n+2, \quad&
b_3&=\beta_3n+2,
\end{alignedat}
\label{eq:6.7}
\end{equation}
and consider, as in Section~\ref{sec:5},
the corresponding set of parameters
\begin{align*}
c_{00}&=(\beta_1+\beta_2+\beta_3)-(\alpha_1+\alpha_2+\alpha_3),
\\
c_{jk}&=\begin{cases}
\alpha_j-\beta_k & \mbox{if $\alpha_j>\beta_k$}, \\
\beta_k-\alpha_j & \mbox{if $\alpha_j<\beta_k$},
\end{cases}
\qquad j,k=1,2,3;
\end{align*}
hence the set $\bc\cdot n$ corresponds to~\eqref{eq:6.7}.
Set
\begin{align*}
m_1=m_1(\bc)
&:=\beta_3^*-\alpha_1^*,
\\
m_2=m_2(\bc)
&:=\max\{\alpha_1-\beta_1,
\beta_3^*-\alpha_2,\beta_3^*-\alpha_3,\beta_2^*-\alpha_1^*\},
\\
m_3=m_3(\bc)
&:=\min\{m_1(\bc),m_2(\bc)\},
\end{align*}
where asterisks mean ordering in accordance with~\eqref{eq:6.4}.
To the $60$-element set
$\sM=\sM(\bc)=\{\fq\,\bc:\fq\in\fG/\fG_0\}$
we assign the function
\begin{align*}
\phi(x)=\smash{\max_{\bc'\in\sM}}&\bigl(
\[c_{00}x\]+\[c_{21}x\]+\[c_{31}x\]+\[c_{22}x\]+\[c_{33}x\]
\\ &\qquad
-\[c_{00}'x\]-\[c_{21}'x\]-\[c_{31}'x\]-\[c_{22}'x\]-\[c_{33}'x\]\bigr),
\end{align*}
which is $1$-periodic and takes only non-negative integral values.
Further, let $\tau_0$ and~$\tau_1$, $\tau_0<\tau_1$, be
the (real) zeros  of the quadratic polynomial
$$
(\tau-\alpha_1)(\tau-\alpha_2)(\tau-\alpha_3)
-(\tau-\beta_1)(\tau-\beta_2)(\tau-\beta_3)
$$
(in particular, $\tau_0<\beta_1$ and $\tau_1>\alpha_3^*$)
and let
\begin{align*}
f_0(\tau)
&=\alpha_1\log(\alpha_1-\tau)
+\alpha_2\log(\alpha_2-\tau)
+\alpha_3\log(\alpha_3-\tau)
\\ &\qquad
-\beta_1\log(\tau-\beta_1)
-\beta_2\log(\beta_2-\tau)
-\beta_3\log(\beta_3-\tau)
\\ &\qquad
-(\alpha_1-\beta_1)\log(\alpha_1-\beta_1)
+(\beta_2-\alpha_2)\log(\beta_2-\alpha_2)
\\ &\qquad
+(\beta_3-\alpha_3)\log(\beta_3-\alpha_3)
\end{align*}
be a function in the cut $\tau$-plane
$\mathbb C\setminus(-\infty,\beta_1]\cup[\alpha_1^*,+\infty)$.
Then the final result is as follows.

\begin{proposition}
\label{prop:4}
In the above notation let
$$
\begin{gathered}
C_0=-\Re f_0(\tau_0), \qquad C_1=\Re f_0(\tau_1),
\\
C_2=m_1+m_2
-\biggl(\int_0^1\phi(x)\,\d\psi(x)
-\int_0^{1/m_3}\phi(x)\frac{\d x}{x^2}\biggr).
\end{gathered}
$$
If $C_0>C_2$\rom, then
$$
\mu(\zeta(2))\le\frac{C_0+C_1}{C_0-C_2}.
$$
\end{proposition}

In accordance with~\cite{RV2} we now take
\begin{equation}
\begin{alignedat}{3}{3}
\alpha_1&=13, \quad& \alpha_2&=12, \quad& \alpha_3&=14,
\\
\beta_1&=0, \quad& \beta_2&=24, \quad& \beta_3&=28
\end{alignedat}
\label{eq:6.8}
\end{equation}
and obtain the following result.

\begin{theorem}[\textrm{\rm\cite{RV2}}]
\label{th:2}
The irrationality exponent of~$\zeta(2)$ satisfies the estimate
\begin{equation}
\mu(\zeta(2))\le5.44124250\dots\,.
\label{eq:6.9}
\end{equation}
\end{theorem}

\begin{observation}
In addition to the fact that the group for~$\zeta(2)$ can be
naturally embedded into the group for~$\zeta(3)$, we can make the following
surprising observation relating the best known estimates of the
irrationality exponents for these constants. The choice
of the directions~\eqref{eq:5.1} with
$$
\begin{alignedat}{4}{4}
\alpha_1&=16, \quad& \alpha_2&=17, \quad&
\alpha_3&=18, \quad& \alpha_4&=19,
\\
\beta_1&=0, \quad& \beta_2&=7, \quad&
\beta_3&=31, \quad& \beta_4&=32
\end{alignedat}
$$
for~$\zeta(3)$ (cf.~\eqref{eq:5.12}\,)
and the choice of the directions~\eqref{eq:6.7} with
$$
\begin{alignedat}{3}{3}
\alpha_1&=10, \quad& \alpha_2&=11, \quad& \alpha_3&=12,
\\
\beta_1&=0, \quad& \beta_2&=24, \quad& \beta_3&=25
\end{alignedat}
$$
for~$\zeta(2)$
(which is $\fG$-equivalent to~\eqref{eq:6.8}\,)
lead to the following matrices~\eqref{eq:4.8} and~\eqref{eq:6.6}:
\begin{equation}
\pmatrix
16 &  9 & 15 & 16 \\
17 & 10 & 14 & 15 \\
18 & 11 & 13 & 14 \\
19 & 12 & 12 & 13
\endpmatrix
\qquad\mbox{and}\qquad
\pmatrix
16 &    &    &    \\
   & 10 & 14 & 15 \\
   & 11 & 13 & 14 \\
   & 12 & 12 & 13
\endpmatrix.
\label{eq:6.10}
\end{equation}
The first set of the parameters in~\eqref{eq:6.10} produces
the estimate~\eqref{eq:5.11}, while the second set the
estimate~\eqref{eq:6.9}.
\end{observation}

Finally, we point out that the known group structure for~$\log2$
(and for some other values of the Gauss hypergeometric function)
is quite simple since no identity like~\eqref{eq:4.1} is known;
the corresponding group consists of just two permutations
(see~\cite{Vi} for an explanation in terms of `multiple' integrals).

\section{Arithmetic of special rational functions}
\label{sec:7}
In our study of arithmetic properties of linear forms in~$1$ and~$\zeta(3)$
we have used the information coming mostly from
$G$-presentations~\eqref{eq:4.13}.
If we denote by~$F(\bh)$ the right-hand side of~\eqref{eq:5.9}
and apply Lemma~\ref{lem:7}, then one could think that the expansion
\begin{equation}
F(\bh)=\sum_{t=0}^\infty R(t),
\label{eq:7.1}
\end{equation}
where we now set
$$
R(t)=R(h_0;h_1,h_2,h_3,h_4,h_5;t)
=(h_0+2t)\prod_{j=1}^6R_j(t)
$$
with
\begin{equation}
\begin{gathered}
\begin{aligned}
R_1(t)&=(h_0-h_1-h_2)!\cdot\frac{\Gamma(h_1+t)}{\Gamma(1+h_0-h_2+t)},
\\
R_2(t)&=(h_0-h_2-h_4)!\cdot\frac{\Gamma(h_2+t)}{\Gamma(1+h_0-h_4+t)},
\\
R_3(t)&=(h_0-h_1-h_3)!\cdot\frac{\Gamma(h_3+t)}{\Gamma(1+h_0-h_1+t)},
\\
R_4(t)&=(h_0-h_3-h_5)!\cdot\frac{\Gamma(h_5+t)}{\Gamma(1+h_0-h_3+t)},
\end{aligned}
\\
R_5(t)=\frac1{(h_4-1)!}\cdot\frac{\Gamma(h_4+t)}{\Gamma(1+t)},
\quad
R_6(t)=\frac1{(h_5-1)!}\cdot\frac{\Gamma(h_0+t)}{\Gamma(1+h_0-h_5+t)},
\end{gathered}
\label{eq:7.2}
\end{equation}
brings with it some extra arithmetic for linear forms~$H(\bc)$
since the functions~\eqref{eq:7.2} are of the same type as~\eqref{eq:2.5}.
Unfortunately, we have discovered that (quite complicated from the
computational point of view) arithmetic of the presentations~\eqref{eq:7.1}
brings nothing new.

For our future aims we now study the arithmetic properties
of elementary `bricks'---rational functions
\begin{equation}
R(t)=R(a,b;t):=\begin{cases}
\dfrac{(t+b)(t+b+1)\dotsb(t+a-1)}{(a-b)!}
& \mbox{if $a\ge b$}, \\
\dfrac{(b-a-1)!}{(t+a)(t+a+1)\dotsb(t+b-1)}
& \mbox{if $a<b$},
\end{cases}
\label{eq:7.3}
\end{equation}
which are introduced by Nesterenko~\cite{Ne2,Ne3}
and appear in~\eqref{eq:2.5} and~\eqref{eq:7.2}.

The next claim exploits well-known properties of integral-valued
polynomials.

\begin{lemma}[\textrm{\rm cf.~Lemma~\ref{lem:1}}]
\label{lem:15}
Suppose that $a\ge b$. Then for any non-negative integer~$j$
there hold the inclusions
$$
D_{a-b}^j\cdot\frac1{j!}R^{(j)}(-k)\in\mathbb Z, \qquad k\in\mathbb Z.
$$
\end{lemma}

The next claim immediately follows from Lemma~\ref{lem:2}
in the same way as Lemma~\ref{lem:3}.

\begin{lemma}
\label{lem:16}
Let $a,b,a_0,b_0$ be integers, $a_0\le a<b\le b_0$.
Then for any non-negative integer~$j$
there hold the inclusions
$$
D_{b_0-a_0-1}^j\cdot\frac1{j!}\bigl(R(t)(t+k)\bigr)^{(j)}\big|_{t=-k}
\in\mathbb Z,
\qquad k=a_0,a_0+1,\dots,b_0-1.
$$
\end{lemma}

Lemmas \ref{lem:15} and \ref{lem:16}
give a particular (but quite important)
information on the $p$-adic valuation of the values
$R^{(j)}(-k)$ and $\bigl(R(t)(t+k)\bigr)^{(j)}\big|_{t=-k}$
respectively, with a help of the formula
$\ord_pD_N=1$ for any integer~$N$ and any prime~$p$
in the interval $\sqrt N<p\le N$.
Two next statements are devoted to
the `most precise' estimates for the $p$-adic order
of these quantities.

\begin{lemma}
\label{lem:17}
Let $a,b,a_0,b_0$~be integers\rom, $b_0\le b<a\le a_0$\rom, and
let $R(t)=R(a,b;t)$~be defined by~\eqref{eq:7.3}.
Then for any integer~$k$\rom, $b_0\le k<a_0$\rom,
any prime $p>\sqrt{a_0-b_0-1}$\rom, and any non-negative integer~$j$
there hold the estimates
\begin{align}
\ord_pR^{(j)}(-k)
&\ge-j+\biggl\[\frac{a-1-k}p\biggr\]
-\biggl\[\frac{b-1-k}p\biggr\]-\biggl\[\frac{a-b}p\biggr\]
\nonumber\\
&=-j+\biggl\[\frac{k-b}p\biggr\]
-\biggl\[\frac{k-a}p\biggr\]-\biggl\[\frac{a-b}p\biggr\].
\label{eq:7.4}
\end{align}
\end{lemma}

\begin{proof}
Fix an arbitrary prime $p>\sqrt{a_0-b_0-1}$.
First, we note that by the definition of the integral part of a number
$$
\[-x\]=-\[x\]-\delta_x,
\qquad\mbox{where}\quad
\delta_x=\begin{cases}
0 & \mbox{if $x\in\mathbb Z$}, \\
1 & \mbox{if $x\notin\mathbb Z$},
\end{cases}
$$
which yields
$$
\biggl\[-\frac sp\biggr\]=-\biggl\[\frac{s-1}p\biggr\]-1
\qquad\mbox{for}\quad s\in\mathbb Z.
$$
Therefore,
\begin{equation}
\biggl\[\frac{k-b}p\biggr\]
=-\biggl\[\frac{b-1-k}p\biggr\]-1,
\qquad
\biggl\[\frac{a-1-k}p\biggr\]
=-\biggl\[\frac{k-a}p\biggr\]-1
\label{eq:7.5}
\end{equation}
for any integer~$k$.

Direct calculations show that
$$
R(-k)=\begin{cases}
\dfrac{(a-1-k)!}{(b-1-k)!\,(a-b)!}
& \mbox{if $k<b$}, \\
0
& \mbox{if $b\le k<a$}, \\
(-1)^{a-b}\dfrac{(k-b)!}{(k-a)!\,(a-b)!}
& \mbox{if $k\ge a$};
\end{cases}
$$
thus,
$$
\begin{alignedat}{2}{2}
\ord_pR(-k)
&\ge\biggl\[\frac{a-1-k}p\biggr\]
-\biggl\[\frac{b-1-k}p\biggr\]-\biggl\[\frac{a-b}p\biggr\]
\qquad&& \mbox{if $k<a$},
\\
\ord_pR(-k)
&\ge\biggl\[\frac{k-b}p\biggr\]
-\biggl\[\frac{k-a}p\biggr\]-\biggl\[\frac{a-b}p\biggr\]
\qquad&& \mbox{if $k\ge b$},
\end{alignedat}
$$
which yields the estimates~\eqref{eq:7.4} for $j=0$ with the help
of~\eqref{eq:7.5}.

If $k<b$ or $k\ge a$, consider the function
$$
r(t)=\frac{R'(t)}{R(t)}
=\sum_{l=b}^{a-1}\frac1{t+l},
$$
hence for any integer $j\ge1$ there hold the inclusions
$$
r^{(j-1)}(-k)\cdot D_{\max\{a-b_0-1,a_0-b-1\}}^{j-1}\in\mathbb Z.
$$
Induction on~$j$ and the identity
\begin{equation}
R^{(j)}(t)
=\bigl(R(t)r(t)\bigr)^{(j-1)}
=\sum_{m=0}^{j-1}\binom{j-1}mR^{(m)}(t)r^{(j-1-m)}(t)
\label{eq:7.6}
\end{equation}
specified at $t=-k$ lead us to the required estimates~\eqref{eq:7.4}.

If $b\le k<a$, consider the functions
$$
R_k(t)=\frac{R(t)}{t+k},
\qquad
r_k(t)=\frac{R_k'(t)}{R_k(t)}
=\sum\doublesb{l=b}{l\ne k}^{a-1}\frac1{t+l};
$$
obviously, for any integer $j\ge1$ there hold the inclusions
$$
r_k^{(j-1)}(-k)\cdot D_{a-b-1}^{j-1}\in\mathbb Z.
$$
Then
$$
R^{(j)}(-k)=jR_k^{(j-1)}(-k)
$$
since
$$
R_k(-k)=(-1)^{k-b}\frac{(k-b)!\,(a-1-k)!}{(a-b)!},
$$
and induction on~$j$ in combination with identity~\eqref{eq:7.6}
(where we substitute $R_k(t),\linebreak[2]r_k(t)$
for $R(t),r(t)$, respectively)
show that
\begin{align*}
\ord_pR^{(j)}(-k)
&\ge\ord_pR_k^{(j-1)}(-k)
\\
&\ge-(j-1)+\biggl\[\frac{k-b}p\biggr\]
+\biggl\[\frac{a-1-k}p\biggr\]-\biggl\[\frac{a-b}p\biggr\]
\end{align*}
for integer $j\ge1$. Thus, applying~\eqref{eq:7.5}
we obtain the required estimates~\eqref{eq:7.4} again.
The proof is complete.
\end{proof}

\begin{lemma}
\label{lem:18}
Let $a,b,a_0,b_0$~be integers\rom, $a_0\le a<b\le b_0$\rom, and
let $R(t)= R(a,b;t)$~be defined by~\eqref{eq:7.3}.
Then for any integer~$k$\rom, $a_0\le k<b_0$\rom,
any prime $p>\sqrt{b_0-a_0-1}$\rom, and any non-negative integer~$j$
there hold the estimates
\begin{equation}
\ord_p\bigl(R(t)(t+k)\bigr)^{(j)}\big|_{t=-k}
\ge-j+\biggl\[\frac{b-a-1}p\biggr\]
-\biggl\[\frac{k-a}p\biggr\]-\biggl\[\frac{b-1-k}p\biggr\].
\label{eq:7.7}
\end{equation}
\end{lemma}

\begin{proof}
Fix an arbitrary prime $p>\sqrt{b_0-a_0-1}$. We have
$$
\bigl(R(t)(t+k)\bigr)\big|_{t=-k}=\begin{cases}
(-1)^{k-a}\dfrac{(b-a-1)!}{(k-a)!\,(b-1-k)!}
& \mbox{if $a\le k<b$}, \\
0 & \mbox{if $k<a$ or $k\ge b$},
\end{cases}
$$
which yields the estimates~\eqref{eq:7.7} for $j=0$.

Considering in the case $a\le k<b$ the functions
$$
R_k(t)=R(t)(t+k),
\qquad
r_k(t)=\frac{R_k'(t)}{R_k(t)}
=\sum\doublesb{l=a}{l\ne k}^{b-1}\frac1{t+l},
$$
and carrying out induction on $j\ge0$,
with the help of identity~\eqref{eq:7.6}
(where we take $R_k(t),r_k(t)$ for
$R(t),r(t)$ again) we deduce the estimates~\eqref{eq:7.7}.

If $k<a$ or $k\ge b$ note that
$$
\bigl(R(t)(t+k)\bigr)^{(j)}\big|_{t=-k}
=jR^{(j-1)}(-k).
$$
Since
$$
R(-k)=\begin{cases}
\dfrac{(b-a-1)!\,(a-1-k)!}{(b-1-k)!}
& \mbox{if $k<a$}, \\
(-1)^{b-a}\dfrac{(b-a-1)!\,(k-b)!}{(k-a)!}
& \mbox{if $k\ge b$},
\end{cases}
$$
induction on~$j$ and equalities~\eqref{eq:7.5} yield
the required estimates~\eqref{eq:7.7} again.
The proof is complete.
\end{proof}

\section{Linear forms in~$1$ and odd zeta values}
\label{sec:8}
Since generalizations of $G$-presentations~\eqref{eq:2.13},~\eqref{eq:6.3}
lead us to forms involving both odd and even zeta values, it is
natural to follow Rivoal dealing with $F$-presentations.

Consider positive odd integers $q$ and $r$, where $q\ge r+4$.
To a set of integral positive parameters
$$
\bh=(h_0;h_1,\dots,h_q)
$$
satisfying the condition
\begin{equation}
h_1+h_2+\dots+h_q\le h_0\cdot\frac{q-r}2
\label{eq:8.1}
\end{equation}
we assign the rational function
\begin{align}
\wt R(t)
&\phantom:=\wt R(\bh;t)
\nonumber\\
&:=(h_0+2t)\frac{\Gamma(h_0+t)^r\Gamma(h_1+t)\dotsb\Gamma(h_q+t)}
{\Gamma(1+t)^r\Gamma(1+h_0-h_1+t)\dotsb\Gamma(1+h_0-h_q+t)}.
\label{eq:8.2}
\end{align}
By~\eqref{eq:8.1} we obtain
\begin{equation}
\wt R(t)=O\biggl(\frac1{t^2}\biggr),
\label{eq:8.3}
\end{equation}
hence the quantity
\begin{equation}
\wt F(\bh)
:=\frac1{(r-1)!}\sum_{t=0}^\infty\wt R^{(r-1)}(t)
\label{eq:8.4}
\end{equation}
is well-defined. If $r=1$, the quantity~\eqref{eq:8.4} can be written
as a well-poised hypergeometric series with a special form of the second
parameter; namely,
\begin{align*}
\wt F(\bh)
&=\frac{h_0!\,(h_1-1)!\dotsb(h_q-1)!}
{(h_0-h_1)!\dotsb(h_0-h_q)!}
\\ &\qquad\times 
{}_{q+2}\!F_{q+1}\biggl(\begin{array}{rrrrr}
h_0, & 1+\frac12h_0, &       h_1, & \dots, &       h_q \\[1pt]
     &   \frac12h_0, & 1+h_0-h_1, & \dots, & 1+h_0-h_q
\end{array}\biggm|1\biggr)
\end{align*}
(cf.~\eqref{eq:4.2}\,), while in the case $r>1$ we obtain
a linear combination of well-poised Meijer's $G$-functions
taken at the points~$e^{\pi ik}$,
where $k=\pm1,\pm3,\dots,\pm(r-2)$.

Applying the symmetry of the rational function~\eqref{eq:8.2}
under the substitution $t\mapsto-t-h_0$:
\begin{equation}
\wt R(-t-h_0)=-(-1)^{h_0(q+r)}\wt R(t)=-\wt R(t),
\label{eq:8.5}
\end{equation}
where we use the identity~\eqref{eq:3.4},
and following the arguments of the proof of Lemma~\ref{lem:4}
we are now able to state that the quantity~\eqref{eq:8.4}
is a linear form in~$1$ and odd zeta values
with rational coefficients.
To present this result explicitly we require the ordering
$$
h_1\le h_2\le\dots\le h_q<\frac12h_0
$$
and the following arithmetic normalization of~\eqref{eq:8.4}:
\begin{equation}
F(\bh)
:=\frac{\prod_{j=r+1}^q(h_0-2h_j)!}
{\prod_{j=1}^r(h_j-1)!^2}\cdot\wt F(\bh)
=\frac1{(r-1)!}\sum_{t=1-h_1}^\infty R^{(r-1)}(t),
\label{eq:8.6}
\end{equation}
where the rational function
\begin{equation}
\begin{split}
R(t)
&:=\prod_{j=1}^r\frac1{(h_j-1)!}\,
\frac{\Gamma(h_j+t)}{\Gamma(1+t)}
\cdot\prod_{j=1}^r\frac1{(h_j-1)!}\,
\frac{\Gamma(h_0+t)}{\Gamma(1+h_0-h_j+t)}
\\ &\qquad\times
\prod_{j=r+1}^q(h_0-2h_j)!\,
\frac{\Gamma(h_j+t)}{\Gamma(1+h_0-h_j+t)}
\end{split}
\label{eq:8.7}
\end{equation}
is the product of elementary bricks~\eqref{eq:7.3}.
Set $m_0=\max\{h_r-1,h_0-2h_{r+1}\}$
and $m_j=\max\{m_0,h_0-h_1-h_{r+j}\}$
for $j=1,\dots,q-r$,
and define the integral quantity
\begin{equation}
\Phi=\Phi(\bh):=\prod_{\sqrt{h_0}<p\le m_{q-r}}p^{\nu_p},
\label{eq:8.8}
\end{equation}
where
\begin{equation}
\nu_p
:=\min_{h_{r+1}\le k\le h_0-h_{r+1}}\{\nu_{k,p}\}
\label{eq:8.9}
\end{equation}
and
\begin{align*}
\nu_{k,p}
&:=\sum_{j=1}^r\biggl(\biggl\[\frac{k-1}p\biggr\]
+\biggl\[\frac{h_0-k-1}p\biggr\]
\\ &\phantom:\qquad
-\biggl\[\frac{k-h_j}p\biggr\]
-\biggl\[\frac{h_0-h_j-k}p\biggr\]
-2\biggl\[\frac{h_j-1}p\biggr\]\biggr)
\\ &\phantom:\qquad
+\sum_{j=r+1}^q\biggl(\biggl\[\frac{h_0-2h_j}p\biggr\]
-\biggl\[\frac{k-h_j}p\biggr\]
-\biggl\[\frac{h_0-h_j-k}p\biggr\]\biggr).
\end{align*}
In this notation the result reads as follows.

\begin{lemma}
\label{lem:19}
The quantity~\eqref{eq:8.6} is a linear form in
$1,\zeta(r+2),\zeta(r+4),\dots,\linebreak[4]\zeta(q-4),\zeta(q-2)$
with rational coefficients\rom;
moreover\rom,
$$
D_{m_1}^rD_{m_2}\dotsb D_{m_{q-r}}
\cdot\Phi^{-1}\cdot F(\bh)
\in\mathbb Z\zeta(q-2)+\mathbb Z\zeta(q-4)+\dots
+\mathbb Z\zeta(r+2)+\mathbb Z.
$$
\end{lemma}

\begin{proof}
Applying the Leibniz rule for differentiating a product,
Lemmas~\ref{lem:15}, \ref{lem:16}
and Lemmas~\ref{lem:17},~\ref{lem:18}
to the rational function~\eqref{eq:8.7}
we see that the numbers
$$
\begin{gathered}
B_{jk}=\frac1{(q-j)!}
\cdot\bigl(R(t)(t+k)^{q-r}\bigr)^{(q-j)}\big|_{t=-k},
\\
j=r+1,\dots,q, \quad k=h_{r+1},\dots,h_0-h_{r+1},
\end{gathered}
$$
satisfy the relations
\begin{equation}
D_{m_0}^{q-j}\cdot B_{jk}\in\mathbb Z
\label{eq:8.10}
\end{equation}
and
\begin{equation}
\ord_pB_{jk}
\ge-(q-j)+\nu_{k,p},
\label{eq:8.11}
\end{equation}
respectively, for any $k=h_{r+1},\dots,h_0-h_{r+1}$
and any prime~$p>\sqrt{h_0}$.
Furthermore, the expansion
$$
R(t)=\sum_{j=r+1}^q\sum_{k=h_j}^{h_0-h_j}\frac{B_{jk}}{(t+k)^{j-r}}
$$
leads us to the series
\begin{align*}
F(\bh)
&=\sum_{j=r+1}^q\binom{j-2}{r-1}
\sum_{k=h_j}^{h_0-h_j}B_{jk}
\biggl(\sum_{l=1}^\infty-\sum_{l=1}^{k-h_1}\biggr)\frac1{l^{j-1}}
\\
&=\sum_{j=r+1}^qA_{j-1}\zeta(j-1)-A_0,
\end{align*}
where
\begin{align}
A_{j-1}&=\binom{j-2}{r-1}\sum_{k=h_j}^{h_0-h_j}B_{jk},
\qquad j=r+1,\dots,q,
\label{eq:8.12}
\\
A_0&=\sum_{j=r+1}^q\binom{j-2}{r-1}
\sum_{k=h_j}^{h_0-h_j}B_{jk}
\sum_{l=1}^{k-h_1}\frac1{l^{j-1}}.
\nonumber
\end{align}
By~\eqref{eq:8.10} and the inclusions
$$
D_{m_1}^rD_{m_2}\dotsb D_{m_{j-r}}
\cdot\sum_{l=1}^{k-h_1}\frac1{l^{j-1}}\in\mathbb Z
$$
for any $k=h_j,\dots,h_0-h_j$, $j=r+1,\dots,q$,
we obtain the `fairly rough' inclusions
\begin{gather*}
D_{m_0}^{q-j-1}\cdot A_j\in\mathbb Z
\qquad\mbox{for}\quad j=r,r+1,\dots,q-1,
\\
D_{m_1}^rD_{m_2}\dotsb D_{m_{q-r}}
\cdot A_0\in\mathbb Z,
\end{gather*}
which are (in a sense) refined by the estimates~\eqref{eq:8.11}:
$$
\ord_pA_j\ge-(q-j-1)+\nu_p
\qquad\mbox{for $j=0$ and $j=r,r+1,\dots,q-1$}
$$
with exponents~$\nu_p$ defined in~\eqref{eq:8.9}.
To complete the proof we must show that
$$
A_r=0 \qquad\mbox{and}\qquad
A_{r+1}=A_{r+3}=\dots=A_{q-3}=A_{q-1}=0.
$$
The first equality follows from~\eqref{eq:8.3};
by~\eqref{eq:8.5} we obtain
$$
B_{jk}=(-1)^jB_{j,h_0-k}
\qquad\mbox{for $j=r+1,\dots,q$},
$$
which yields $A_{j-1}=0$ for odd~$j$ according to~\eqref{eq:8.12}.
The proof is complete.
\end{proof}

To evaluate the growth of the linear forms~\eqref{eq:8.6}
so constructed we define the set of integral directions
$\Beta=(\eta_0;\eta_1,\dots,\eta_q)$ and the increasing
integral parameter~$n$
related with the parameters~$\bh$ by the formulae
\begin{equation}
h_0=\eta_0n+2 \qquad\mbox{and}\qquad
h_j=\eta_jn+1 \quad\mbox{for $j=1,\dots,q$}.
\label{eq:8.13}
\end{equation}

Consider the auxiliary function
\begin{align*}
f_0(\tau)
&=r\eta_0\log(\eta_0-\tau)
+\sum_{j=1}^q\bigl(\eta_j\log(\tau-\eta_j)
-(\eta_0-\eta_j)\log(\tau-\eta_0+\eta_j)\bigr)
\\ &\qquad
-2\sum_{j=1}^r\eta_j\log\eta_j
+\sum_{j=r+1}^q(\eta_0-2\eta_j)\log(\eta_0-2\eta_j)
\end{align*}
defined in the cut $\tau$-plane
$\mathbb C\setminus(-\infty,\eta_0-\eta_1]\cup[\eta_0,+\infty)$.
The next assertion is deduced by an application of
the saddle-point method and the use of the asymtotics
of the gamma factors in~\eqref{eq:8.7}
(see, e.g., \cite{Zu3}, Section~2, or~\cite{Ri4}).
We underline that no approach in terms
of real multiple integrals is known in the case $r\ge3$.

\begin{lemma}
\label{lem:20}
Let $r=3$ and let $\tau_0$~be a zero of the polynomial
$$
(\tau-\eta_0)^r(\tau-\eta_1)\dotsb(\tau-\eta_q)
-\tau^r(\tau-\eta_0+\eta_1)\dotsb(\tau-\eta_0+\eta_q)
$$
with $\Im\tau_0>0$ and the maximum possible value of~$\Re\tau_0$.
Suppose that $\Re\tau_0<\eta_0$ and
$\Im f_0(\tau_0)\notin\pi\mathbb Z$. Then
$$
\limsup_{n\to\infty}\frac{\log|F(\bh)|}n=\Re f_0(\tau_0).
$$
\end{lemma}

We now take
$$
m_j=\max\{\eta_r,\eta_0-2\eta_{r+1},\eta_0-\eta_1-\eta_{r+j}\}
\qquad\mbox{for}\quad j=1,\dots,q-r
$$
(hence we scale down with factor~$n$ the old parameters).
The asymptotics of the quantity~\eqref{eq:8.8} as $n\to\infty$
can be calculated with the use of the integral-valued
function
\begin{align*}
\phi_0(x,y)
&:=\sum_{j=1}^r\bigl(\[y\]+\[\eta_0x-y\]
-\[y-\eta_jx\]-\[(\eta_0-\eta_j)x-y\]-2\[\eta_jx\]\bigr)
\\ &\qquad
+\sum_{j=r+1}^q\bigl(\[(\eta_0-2\eta_j)x\]
-\[y-\eta_jx\]-\[(\eta_0-\eta_j)x-y\]\bigr),
\end{align*}
which is $1$-periodic with respect to each variable~$x$ and~$y$.
Then by~\eqref{eq:8.9} and~\eqref{eq:8.13} we obtain
$$
\nu_p=\min_{\eta_4n\le k-1\le(\eta_0-\eta_4)n}
\phi_0\biggl(\frac np,\frac{k-1}p\biggr)
\ge\phi\biggl(\frac np\biggr),
$$
where
$$
\phi(x):=\min_{y\in\mathbb R}\phi_0(x,y)
=\min_{0\le y<1}\phi_0(x,y).
$$
Therefore, the final result is as follows.

\begin{proposition}
\label{prop:5}
In the above notation let $r=3$ and
$$
\begin{gathered}
C_0=-\Re f_0(\tau_0),
\\
C_2=rm_1+m_2+\dots+m_{q-r}
-\biggl(\int_0^1\phi(x)\,\d\psi(x)
-\int_0^{1/m_{q-r}}\phi(x)\,\frac{\d x}{x^2}\biggr).
\end{gathered}
$$
If $C_0>C_2$\rom, then at least one of the numbers
$$
\zeta(5), \; \zeta(7), \; \dots, \; \zeta(q-4), \;
\mbox{and\/} \; \zeta(q-2)
$$
is irrational.
\end{proposition}

We are now ready to state the following new result.

\begin{theorem}
\label{th:3}
At least one of the four numbers
$$
\zeta(5), \; \zeta(7), \; \zeta(9), \; \mbox{and\/} \; \zeta(11)
$$
is irrational.
\end{theorem}

\begin{proof}
Taking $r=3$, $q=13$,
$$
\eta_0=91, \qquad \eta_1=\eta_2=\eta_3=27,
\qquad \eta_j=25+j \quad\mbox{for}\; j=4,5,\dots,13,
$$
we obtain $\tau_0=87.47900541\hdots+i\,3.32820690\dots$,
\begin{align*}
C_0&=-\Re f_0(\tau_0)=227.58019641\dots,
\\
C_2&=3\cdot35+34+8\cdot33
-\biggl(\int_0^1\phi(x)\,\d\psi(x)
-\int_0^{1/33}\phi(x)\,\frac{\d x}{x^2}\biggr)
\\
&=226.24944266\dots
\end{align*}
since in this case
$$
\phi(x)=\nu \quad\mbox{if $x\in\Omega_\nu\setminus\Omega_{\nu+1}$},
\qquad \nu=0,1,\dots,9,
$$
for $x\in[0,1)$, where $\Omega_0=[0,1)$,
\begin{align*}
\Omega_1
&=\Omega_2
=\bigl[\tfrac2{91},\tfrac{36}{37}\bigr)
\cup\bigl[\tfrac{90}{91},1\bigr),
\displaybreak[0]\\
\Omega_3
&=\bigl[\tfrac2{91},\tfrac1{20}\bigr)
\cup\bigl[\tfrac5{91},\tfrac34\bigr)
\cup\bigl[\tfrac{28}{37},\tfrac{13}{14}\bigr)
\cup\bigl[\tfrac{14}{15},\tfrac{35}{37}\bigr)
\cup\bigl[\tfrac{18}{19},\tfrac{27}{28}\bigr)
\cup\bigl[\tfrac{88}{91},\tfrac{36}{37}\bigr)
\cup\bigl[\tfrac{90}{91},1\bigr),
\displaybreak[0]\\
\Omega_4
&=\bigl[\tfrac1{38},\tfrac1{22}\bigr)
\cup\bigl[\tfrac5{91},\tfrac3{26}\bigr)
\cup\bigl[\tfrac2{17},\tfrac18\bigr)
\cup\bigl[\tfrac4{31},\tfrac4{27}\bigr)
\cup\bigl[\tfrac5{33},\tfrac7{30}\bigr)
\cup\bigl[\tfrac4{17},\tfrac{12}{37}\bigr)
\cup\bigl[\tfrac{30}{91},\tfrac13\bigr)
\\ &\;\;
\cup\bigl[\tfrac{31}{91},\tfrac38\bigr)
\cup\bigl[\tfrac{14}{37},\tfrac{11}{28}\bigr)
\cup\bigl[\tfrac{13}{33},\tfrac9{22}\bigr)
\cup\bigl[\tfrac7{17},\tfrac{13}{28}\bigr)
\cup\bigl[\tfrac8{17},\tfrac12\bigr)
\cup\bigl[\tfrac{19}{37},\tfrac9{14}\bigr)
\cup\bigl[\tfrac{20}{31},\tfrac23\bigr)
\\ &\;\;
\cup\bigl[\tfrac{21}{31},\tfrac34\bigr)
\cup\bigl[\tfrac{25}{33},\tfrac{11}{14}\bigr)
\cup\bigl[\tfrac{26}{33},\tfrac{23}{28}\bigr)
\cup\bigl[\tfrac{14}{17},\tfrac{23}{27}\bigr)
\cup\bigl[\tfrac{31}{36},\tfrac{25}{27}\bigr)
\cup\bigl[\tfrac{85}{91},\tfrac{35}{37}\bigr)
\cup\bigl[\tfrac{20}{21},\tfrac{26}{27}\bigr)
\\ &\;\;
\cup\bigl[\tfrac{32}{33},\tfrac{34}{35}\bigr),
\displaybreak[0]\\
\Omega_5
&=\bigl[\tfrac1{37},\tfrac1{27}\bigr)
\cup\bigl[\tfrac1{25},\tfrac1{24}\bigr)
\cup\bigl[\tfrac5{91},\tfrac1{18}\bigr)
\cup\bigl[\tfrac2{35},\tfrac2{27}\bigr)
\cup\bigl[\tfrac3{38},\tfrac1{12}\bigr)
\cup\bigl[\tfrac8{91},\tfrac3{34}\bigr)
\cup\bigl[\tfrac2{21},\tfrac19\bigr)
\\ &\;\;
\cup\bigl[\tfrac4{33},\tfrac18\bigr)
\cup\bigl[\tfrac5{38},\tfrac4{27}\bigr)
\cup\bigl[\tfrac3{19},\tfrac16\bigr)
\cup\bigl[\tfrac5{29},\tfrac5{27}\bigr)
\cup\bigl[\tfrac4{21},\tfrac5{26}\bigr)
\cup\bigl[\tfrac6{29},\tfrac29\bigr)
\cup\bigl[\tfrac5{21},\tfrac7{27}\bigr)
\\ &\;\;
\cup\bigl[\tfrac4{15},\tfrac{10}{37}\bigr)
\cup\bigl[\tfrac27,\tfrac3{10}\bigr)
\cup\bigl[\tfrac7{23},\tfrac4{13}\bigr)
\cup\bigl[\tfrac6{19},\tfrac{12}{37}\bigr)
\cup\bigl[\tfrac{30}{91},\tfrac13\bigr)
\cup\bigl[\tfrac{10}{29},\tfrac7{20}\bigr)
\cup\bigl[\tfrac{13}{37},\tfrac5{14}\bigr)
\\ &\;\;
\cup\bigl[\tfrac{33}{91},\tfrac38\bigr)
\cup\bigl[\tfrac8{21},\tfrac5{13}\bigr)
\cup\bigl[\tfrac{13}{33},\tfrac{11}{27}\bigr)
\cup\bigl[\tfrac{12}{29},\tfrac5{12}\bigr)
\cup\bigl[\tfrac8{19},\tfrac{11}{26}\bigr)
\cup\bigl[\tfrac{14}{33},\tfrac{13}{30}\bigr)
\cup\bigl[\tfrac{40}{91},\tfrac49\bigr)
\\ &\;\;
\cup\bigl[\tfrac5{11},\tfrac{11}{24}\bigr)
\cup\bigl[\tfrac{17}{37},\tfrac6{13}\bigr)
\cup\bigl[\tfrac{17}{36},\tfrac{13}{27}\bigr)
\cup\bigl[\tfrac{16}{33},\tfrac12\bigr)
\cup\bigl[\tfrac{16}{31},\tfrac{14}{27}\bigr)
\cup\bigl[\tfrac8{15},\tfrac{19}{35}\bigr)
\cup\bigl[\tfrac{17}{31},\tfrac59\bigr)
\\ &\;\;
\cup\bigl[\tfrac{19}{33},\tfrac{15}{26}\bigr)
\cup\bigl[\tfrac{18}{31},\tfrac{16}{27}\bigr)
\cup\bigl[\tfrac{20}{33},\tfrac{17}{28}\bigr)
\cup\bigl[\tfrac{19}{31},\tfrac{17}{27}\bigr)
\cup\bigl[\tfrac{11}{17},\tfrac23\bigr)
\cup\bigl[\tfrac{17}{25},\tfrac{15}{22}\bigr)
\cup\bigl[\tfrac{20}{29},\tfrac{19}{27}\bigr)
\\ &\;\;
\cup\bigl[\tfrac{12}{17},\tfrac{17}{24}\bigr)
\cup\bigl[\tfrac{21}{29},\tfrac{20}{27}\bigr)
\cup\bigl[\tfrac{23}{31},\tfrac34\bigr)
\cup\bigl[\tfrac{69}{91},\tfrac79\bigr)
\cup\bigl[\tfrac{15}{19},\tfrac{19}{24}\bigr)
\cup\bigl[\tfrac45,\tfrac{22}{27}\bigr)
\cup\bigl[\tfrac{14}{17},\tfrac{23}{27}\bigr)
\\ &\;\;
\cup\bigl[\tfrac{25}{29},\tfrac{19}{22}\bigr)
\cup\bigl[\tfrac{27}{31},\tfrac78\bigr)
\cup\bigl[\tfrac{29}{33},\tfrac89\bigr)
\cup\bigl[\tfrac{26}{29},\tfrac9{10}\bigr)
\cup\bigl[\tfrac{28}{31},\tfrac{25}{27}\bigr)
\cup\bigl[\tfrac{31}{33},\tfrac{35}{37}\bigr)
\cup\bigl[\tfrac{87}{91},\tfrac{26}{27}\bigr)
\\ &\;\;
\cup\bigl[\tfrac{32}{33},\tfrac{33}{34}\bigr),
\displaybreak[0]\\
\Omega_6
&=\bigl[\tfrac1{36},\tfrac1{27}\bigr)
\cup\bigl[\tfrac1{17},\tfrac2{27}\bigr)
\cup\bigl[\tfrac9{91},\tfrac4{37}\bigr)
\cup\bigl[\tfrac{10}{91},\tfrac19\bigr)
\cup\bigl[\tfrac{12}{91},\tfrac4{27}\bigr)
\cup\bigl[\tfrac{16}{91},\tfrac5{27}\bigr)
\cup\bigl[\tfrac{19}{91},\tfrac8{37}\bigr)
\\ &\;\;
\cup\bigl[\tfrac5{23},\tfrac29\bigr)
\cup\bigl[\tfrac7{29},\tfrac9{37}\bigr)
\cup\bigl[\tfrac{23}{91},\tfrac7{27}\bigr)
\cup\bigl[\tfrac27,\tfrac8{27}\bigr)
\cup\bigl[\tfrac{29}{91},\tfrac{12}{37}\bigr)
\cup\bigl[\tfrac{30}{91},\tfrac13\bigr)
\cup\bigl[\tfrac{33}{91},\tfrac{10}{27}\bigr)
\\ &\;\;
\cup\bigl[\tfrac{15}{38},\tfrac{11}{27}\bigr)
\cup\bigl[\tfrac37,\tfrac{16}{37}\bigr)
\cup\bigl[\tfrac{40}{91},\tfrac49\bigr)
\cup\bigl[\tfrac9{19},\tfrac{13}{27}\bigr)
\cup\bigl[\tfrac{47}{91},\tfrac{14}{27}\bigr)
\cup\bigl[\tfrac7{13},\tfrac{20}{37}\bigr)
\cup\bigl[\tfrac{50}{91},\tfrac59\bigr)
\\ &\;\;
\cup\bigl[\tfrac{53}{91},\tfrac{16}{27}\bigr)
\cup\bigl[\tfrac8{13},\tfrac{23}{37}\bigr)
\cup\bigl[\tfrac{57}{91},\tfrac{17}{27}\bigr)
\cup\bigl[\tfrac{59}{91},\tfrac{24}{37}\bigr)
\cup\bigl[\tfrac{15}{23},\tfrac{17}{26}\bigr)
\cup\bigl[\tfrac{23}{35},\tfrac23\bigr)
\cup\bigl[\tfrac9{13},\tfrac{26}{37}\bigr)
\\ &\;\;
\cup\bigl[\tfrac{64}{91},\tfrac{19}{27}\bigr)
\cup\bigl[\tfrac{66}{91},\tfrac{19}{26}\bigr)
\cup\bigl[\tfrac{67}{91},\tfrac{20}{27}\bigr)
\cup\bigl[\tfrac{13}{17},\tfrac79\bigr)
\cup\bigl[\tfrac45,\tfrac{22}{27}\bigr)
\cup\bigl[\tfrac{76}{91},\tfrac{31}{37}\bigr)
\cup\bigl[\tfrac{16}{19},\tfrac{23}{27}\bigr)
\\ &\;\;
\cup\bigl[\tfrac{29}{33},\tfrac89\bigr)
\cup\bigl[\tfrac{31}{34},\tfrac{34}{37}\bigr)
\cup\bigl[\tfrac{23}{25},\tfrac{25}{27}\bigr)
\cup\bigl[\tfrac{31}{33},\tfrac{33}{35}\bigr)
\cup\bigl[\tfrac{87}{91},\tfrac{26}{27}\bigr),
\displaybreak[0]\\
\Omega_7
&=\bigl[\tfrac1{33},\tfrac1{27}\bigr)
\cup\bigl[\tfrac1{17},\tfrac2{27}\bigr)
\cup\bigl[\tfrac9{91},\tfrac4{37}\bigr)
\cup\bigl[\tfrac{10}{91},\tfrac19\bigr)
\cup\bigl[\tfrac{12}{91},\tfrac5{37}\bigr)
\cup\bigl[\tfrac17,\tfrac4{27}\bigr)
\cup\bigl[\tfrac{16}{91},\tfrac5{27}\bigr)
\\ &\;\;
\cup\bigl[\tfrac{19}{91},\tfrac8{37}\bigr)
\cup\bigl[\tfrac{20}{91},\tfrac29\bigr)
\cup\bigl[\tfrac{22}{91},\tfrac9{37}\bigr)
\cup\bigl[\tfrac9{35},\tfrac7{27}\bigr)
\cup\bigl[\tfrac27,\tfrac8{27}\bigr)
\cup\bigl[\tfrac{29}{91},\tfrac9{28}\bigr)
\cup\bigl[\tfrac{10}{31},\tfrac{11}{34}\bigr)
\\ &\;\;
\cup\bigl[\tfrac{33}{91},\tfrac{10}{27}\bigr)
\cup\bigl[\tfrac{36}{91},\tfrac{15}{37}\bigr)
\cup\bigl[\tfrac{37}{91},\tfrac{11}{27}\bigr)
\cup\bigl[\tfrac37,\tfrac{16}{37}\bigr)
\cup\bigl[\tfrac{40}{91},\tfrac49\bigr)
\cup\bigl[\tfrac{10}{21},\tfrac{13}{27}\bigr)
\cup\bigl[\tfrac{47}{91},\tfrac{14}{27}\bigr)
\\ &\;\;
\cup\bigl[\tfrac7{13},\tfrac{20}{37}\bigr)
\cup\bigl[\tfrac{50}{91},\tfrac59\bigr)
\cup\bigl[\tfrac{53}{91},\tfrac{16}{27}\bigr)
\cup\bigl[\tfrac8{13},\tfrac{23}{37}\bigr)
\cup\bigl[\tfrac{57}{91},\tfrac{17}{27}\bigr)
\cup\bigl[\tfrac{59}{91},\tfrac{24}{37}\bigr)
\cup\bigl[\tfrac9{13},\tfrac{26}{37}\bigr)
\\ &\;\;
\cup\bigl[\tfrac{64}{91},\tfrac{19}{27}\bigr)
\cup\bigl[\tfrac{66}{91},\tfrac{27}{37}\bigr)
\cup\bigl[\tfrac{67}{91},\tfrac{20}{27}\bigr)
\cup\bigl[\tfrac{10}{13},\tfrac79\bigr)
\cup\bigl[\tfrac{73}{91},\tfrac{30}{37}\bigr)
\cup\bigl[\tfrac{74}{91},\tfrac{22}{27}\bigr)
\cup\bigl[\tfrac{11}{13},\tfrac{23}{27}\bigr)
\\ &\;\;
\cup\bigl[\tfrac{80}{91},\tfrac89\bigr)
\cup\bigl[\tfrac{83}{91},\tfrac{34}{37}\bigr)
\cup\bigl[\tfrac{12}{13},\tfrac{25}{27}\bigr)
\cup\bigl[\tfrac{87}{91},\tfrac{26}{27}\bigr),
\displaybreak[0]\\
\Omega_8
&=\bigl[\tfrac1{31},\tfrac1{27}\bigr)
\cup\bigl[\tfrac6{91},\tfrac2{27}\bigr)
\cup\bigl[\tfrac9{91},\tfrac1{10}\bigr)
\cup\bigl[\tfrac3{29},\tfrac4{37}\bigr)
\cup\bigl[\tfrac{10}{91},\tfrac19\bigr)
\cup\bigl[\tfrac2{15},\tfrac5{37}\bigr)
\cup\bigl[\tfrac17,\tfrac4{27}\bigr)
\\ &\;\;
\cup\bigl[\tfrac3{17},\tfrac5{28}\bigr)
\cup\bigl[\tfrac7{38},\tfrac5{27}\bigr)
\cup\bigl[\tfrac7{33},\tfrac8{37}\bigr)
\cup\bigl[\tfrac{20}{91},\tfrac29\bigr)
\cup\bigl[\tfrac8{33},\tfrac9{37}\bigr)
\cup\bigl[\tfrac9{31},\tfrac7{24}\bigr)
\cup\bigl[\tfrac5{17},\tfrac8{27}\bigr)
\\ &\;\;
\cup\bigl[\tfrac4{11},\tfrac{10}{27}\bigr)
\cup\bigl[\tfrac{37}{91},\tfrac{11}{27}\bigr)
\cup\bigl[\tfrac{11}{23},\tfrac{13}{27}\bigr)
\cup\bigl[\tfrac7{13},\tfrac{20}{37}\bigr)
\cup\bigl[\tfrac{16}{29},\tfrac59\bigr)
\cup\bigl[\tfrac{53}{91},\tfrac7{12}\bigr)
\cup\bigl[\tfrac{17}{29},\tfrac{16}{27}\bigr)
\\ &\;\;
\cup\bigl[\tfrac{13}{21},\tfrac{23}{37}\bigr)
\cup\bigl[\tfrac{23}{33},\tfrac7{10}\bigr)
\cup\bigl[\tfrac{64}{91},\tfrac{19}{27}\bigr)
\cup\bigl[\tfrac{14}{19},\tfrac{20}{27}\bigr)
\cup\bigl[\tfrac{10}{13},\tfrac{27}{35}\bigr)
\cup\bigl[\tfrac{25}{31},\tfrac{30}{37}\bigr)
\cup\bigl[\tfrac{74}{91},\tfrac{22}{27}\bigr)
\\ &\;\;
\cup\bigl[\tfrac{11}{13},\tfrac{23}{27}\bigr)
\cup\bigl[\tfrac{80}{91},\tfrac{31}{35}\bigr)
\cup\bigl[\tfrac{83}{91},\tfrac{11}{12}\bigr)
\cup\bigl[\tfrac{12}{13},\tfrac{25}{27}\bigr)
\cup\bigl[\tfrac{22}{23},\tfrac{26}{27}\bigr),
\displaybreak[0]\\
\Omega_9
&=\bigl[\tfrac1{29},\tfrac1{28}\bigr)
\cup\bigl[\tfrac2{29},\tfrac1{14}\bigr)
\cup\bigl[\tfrac7{19},\tfrac{10}{27}\bigr)
\cup\bigl[\tfrac{12}{25},\tfrac{13}{27}\bigr)
\cup\bigl[\tfrac{17}{23},\tfrac{20}{27}\bigr)
\cup\bigl[\tfrac{15}{17},\tfrac{23}{26}\bigr)
\cup\bigl[\tfrac{24}{25},\tfrac{25}{26}\bigr),
\end{align*}
and $\Omega_{10}=\emptyset$.

The application of Proposition~\ref{prop:5} completes the proof.
\end{proof}

\begin{rem}
In~\cite{Zu4} we consider a particular case
of the above construction and arrive at the irrationality
of at least one of the eight odd zeta values
starting from~$\zeta(5)$;
namely, we take $r=3$, $q=21$, $\eta_0=20$, and
$\eta_1=\dots=\eta_{21}=7$ to achieve this result.
\end{rem}

Looking over all integral directions
$\Beta=(\eta_0;\eta_1,\dots,\eta_q)$ with $q=7$, $9$, and~$11$
satisfying the conditions
$$
\eta_1\le\eta_2\le\dots\le\eta_q<\frac12\eta_0
\qquad\mbox{and}\qquad
\eta_0\le120
$$
we have discovered that no set~$\Beta$ yields the irrationality of
at least one of the numbers $\zeta(5)$, $\zeta(7)$, and~$\zeta(9)$
via Proposition~\ref{prop:5}. Thus, we can think about natural bounds
of the `pure' arithmetic approach achieved in Theorem~\ref{th:3}.

In a similar way our previous results~\cite{Zu4} on the irrationality
of at least one of the numbers in each of the two sets
\begin{gather*}
\zeta(7), \; \zeta(9), \; \zeta(11), \;
\dots, \; \zeta(33), \; \zeta(35),
\\
\zeta(9), \; \zeta(11), \; \zeta(13), \;
\dots, \; \zeta(49), \; \zeta(51)
\end{gather*}
can be improved. We are not able to demonstrate the general
case of Lemma~\ref{lem:20}, although this lemma
(after removing the hypothesis
$\Re\tau_0<\eta_0$) remains true for odd $r>3$ and
for any suitable choice of directions~$\Beta$
(cf.~\cite{Zu3}, Section~2).

\section{One arithmetic conjecture and\\
group structures for odd zeta values}
\label{sec:9}
To expose the arithmetic of linear forms produced
by the quantities~\eqref{eq:8.4} in the general case
we require a certain normalization
by factorials similar to~\eqref{eq:7.1},~\eqref{eq:7.2}, or~\eqref{eq:8.6}.
To this end we introduce a contiguous set of parameters~$\be$:
\begin{equation}
e_{0k}=h_k-1, \;\; 1\le k\le q,
\quad\mbox{and}\quad
e_{jk}=h_0-h_j-h_k, \;\; 1\le j<k\le q,
\label{eq:9.1}
\end{equation}
which plays the same role as the set~$\bc$
in Sections~\ref{sec:4}--\ref{sec:6},
and fix a normalization
$$
F(\bh)=\frac{\varPi_1(\be)}{\varPi_2(\be)}\wt F(\bh),
$$
where $\varPi_1(\be)$~is a product of some $q-r$~factorials of~$e_{jk}$
and $\varPi_2(\be)$~is a product of $2r$~factorials of~$e_{0k'}$
with indices satisfying the condition
$$
\bigcup_{j,k}\{j,k\}\cup\bigcup_{k'}\{k'\}
=\{1,2,\dots,q\}\cup\{1,2,\dots,q\}.
$$
For simplicity we can present a concrete normalization;
denoting
$$
\begin{aligned}
a_j&=\begin{cases}
h_j &\mbox{for $j=1,\dots,q$}, \\
h_0 &\mbox{for $j=q+1,\dots,q+r$},
\end{cases}
\\
b_j&=\begin{cases}
1 &\mbox{for $j=1,\dots,r$}, \\
1+h_0-h_{j-r} &\mbox{for $j=r+1,\dots,r+q$},
\end{cases}
\end{aligned}
$$
we define the rational function
$$
R(t)=R(\bh;t)
:=(h_0+2t)\prod_{j=1}^{q+r}R(a_j,b_j;t)
$$
(where the bricks $R(a_j,b_j;t)$ are defined in~\eqref{eq:7.3}\,)
and the corresponding quantity
\begin{equation}
F(\bh)
:=\frac1{(r-1)!}\sum_{t=0}^\infty R^{(r-1)}(t)
=\frac{\prod_{j=r+1}^qe_{j-r,j}!}
{\prod_{j=1}^re_{0j}!\cdot\prod_{j=q+1}^{q+r}e_{0,j-r}!}
\cdot\wt F(\bh).
\label{eq:9.2}
\end{equation}

Nesterenko's theorem in~\cite{Ne3} (which is not the same as
Proposition~\ref{prop:1} in Section~\ref{sec:3})
and our results in Section~\ref{sec:7}
yield the inclusion
\begin{equation}
D_{m_1}^rD_{m_2}\dotsb D_{m_{q-r}}
\cdot F(\bh)\in\mathbb Z\zeta(q-2)
+\mathbb Z\zeta(q-4)+\dots
+\mathbb Z\zeta(r+2)+\mathbb Z,
\label{eq:9.3}
\end{equation}
where $m_1,m_2,\dots,m_{q-r}$ are the successive maxima
of the set~$\be$, and Lemmas~\ref{lem:17},~\ref{lem:18} allow us to exclude
extra primes appearing in coefficients of linear forms~\eqref{eq:9.3}.

In spite of the natural arithmetic~\eqref{eq:9.3}
of the linear forms~\eqref{eq:9.2}, Ball's example~\eqref{eq:4.3}
supplemented with direct calculations for small
values of~$h_0,h_1,\dots,h_q$
and Rivoal's conjecture~\cite{Ri3}, Section~5.1,
enables us to suggest the following.

\begin{conj}
There holds the inclusion
$$
D_{m_1}^rD_{m_2}\dotsb D_{m_{q-r-1}}
\cdot F(\bh)\in\mathbb Z\zeta(q-2)
+\mathbb Z\zeta(q-4)+\dots
+\mathbb Z\zeta(r+2)+\mathbb Z,
$$
where $m_1,m_2,\dots,m_{q-r-1}$ are the successive maxima
of the set~\eqref{eq:9.1}.
\end{conj}

We underline that a similar conjecture does not hold for the quantities
$$
F(\bh;z):=\frac1{(r-1)!}\sum_{t=0}^\infty R^{(r-1)}(t)z^t
\qquad\mbox{with $z\ne\pm1$}
$$
producing linear forms in polylogarithms;
the case $z=\pm1$ is exceptional.

If this conjecture is true, cancellation of extra primes
with the help of Lemmas~\ref{lem:17},~\ref{lem:18}
becomes almost useless, while the
action of the $\bh$-trivial group (i.e., the group of all permutations
of the parameters $h_1,\dots,h_q$) comes into play. Indeed,
the quantity
$$
\wt F(\bh)=\frac{\varPi_2(\be)}{\varPi_1(\be)}\cdot F(\bh)
$$
is stable under any permutation of $h_1,\dots,h_q$,
hence we can apply arguments similar to the ones
considered in Section~\ref{sec:5}
to cancell extra primes.

Finally, we mention that an analytic evaluation of linear
forms~$F(\bh)$ and their coefficients after a choice of
directions and an increasing parameter~$n$ can be carried out
by the saddle-point method, as in~\cite{Zu3}, Sections~2 and~3
(see also~\cite{He,Ri4,Ne3}).

The particular case $r=1$ of the above construction
can be regarded as a natural generalization of both the Rhin--Viola
approach for~$\zeta(3)$ and Rivoal's construction \cite{Ri1}.
In this case we deal with usual well-poised hypergeometric series,
and the group structure considered above, provided that
Conjecture holds, as well as the approach of Section~\ref{sec:8}
will bring new estimates for the dimensions of the spaces
spanned over~$\mathbb Q$ by~$1$
and~$\zeta(3),\zeta(5),\zeta(7),\dots$\,.
If we set $r=1$, $q=k+2$, $h_0=3n+2$, and $h_1=\dots=h_q=n+1$
in formula~\eqref{eq:9.2}, where $n,k$~are positive
integers and $k\ge3$~is odd, and consider the corresponding
sequence
\begin{equation}
\begin{split}
F_{k,n}
&=2n!^{k-1}\sum_{t=1}^\infty\biggl(t+\frac n2\biggr)
\frac{(t-1)\dotsb(t-n)\cdot(t+n+1)\dotsb(t+2n)}
{t^{k+1}(t+1)^{k+1}\dotsb(t+n)^{k+1}}
\\
&\in\mathbb Q\zeta(k)+\mathbb Q\zeta(k-2)+\dots
+\mathbb Q\zeta(3)+\mathbb Q,
\qquad n=1,2,\dots
\end{split}
\label{eq:9.4}
\end{equation}
(cf.~\eqref{eq:4.3}\,), then it is easy to verify that
\begin{equation}
\lim_{n\to\infty}\frac{\log|F_{5,n}|}n
=-6.38364071\dots\,.
\label{eq:9.5}
\end{equation}
The mysterious thing here is the coincidence of the
asymptotics~\eqref{eq:9.5} of the linear forms~$F_{5,n}$
with the asymptotics of Vasilyev's multiple integrals
$$
J_n(5)=\idotsint\limits_{[0,1]^5}
\frac{x_1^n(1-x_1)^n\dotsb x_5^n(1-x_5)^n
\,\d x_1\dotsb\d x_5}
{(1-(1-(1-(1-(1-x_1)x_2)x_3)x_4)x_5)^{n+1}},
$$
for which the inclusions
$$
D_n^5\cdot J_n(5)\in\mathbb Z\zeta(5)+\mathbb Z\zeta(3)+\mathbb Z,
\qquad n=1,2,\dots,
$$
are proved in~\cite{Va}.
Moreover, we have checked that, numerically,
\begin{gather*}
F_{5,1}=18\zeta(5)+66\zeta(3)-98,
\quad
F_{7,1}=26\zeta(7)+220\zeta(5)+612\zeta(3)-990,
\\
F_{9,1}=34\zeta(9)+494\zeta(7)+2618\zeta(5)+6578\zeta(3)-11154,
\end{gather*}
hence these linear forms are the same forms
as listed in~\cite{Va}, Section~5.
Therefore, it is natural to conjecture%
\footnote{This conjecture is recently proved
in~\cite{Zu6},~\cite{Zu7}.}
the coincidence of Vasilyev's integrals
$$
J_n(k)=\idotsint\limits_{[0,1]^k}
\frac{x_1^n(1-x_1)^nx_2^n(1-x_2)^n\dotsb x_k^n(1-x_k)^n
\,\d x_1\,\d x_2\dotsb\d x_k}
{(1-(1-(\dotsb(1-(1-x_1)x_2)\dotsb)x_{k-1})x_k)^{n+1}},
$$
for odd~$k$
with the corresponding hypergeometric series~\eqref{eq:9.4};
we recall that in the case $k=3$
this coincidence follows from Propositions~\ref{prop:1} and~\ref{prop:2}.
A similar conjecture can be put forward in the case of even~$k$ in view
of Whipple's identity~\eqref{eq:6.5}.

\medskip
We hope that the methods of this work will find a continuation
in the form of new qualitative and quantitative results
on the linear independence of values of the Riemann zeta function
at positive integers.


\end{document}